\documentclass[a4paper]{article}

\usepackage[english]{babel}
\usepackage{amsfonts}
\usepackage[utf8x]{inputenc}
\usepackage[T1]{fontenc}

\usepackage[a4paper,top=3cm,bottom=2cm,left=3cm,right=3cm,marginparwidth=1.75cm]{geometry}

\usepackage{amsmath}
\usepackage{graphicx}
\usepackage{subcaption}
\usepackage[colorinlistoftodos]{todonotes}
\usepackage[colorlinks=true, allcolors=blue]{hyperref}

\title{Energy-preserving Variational Integrators for Forced  Lagrangian Systems}

\author{Harsh Sharma \thanks{Corresponding Author, Email Address for Correspondence: harshs3@vt.edu} \thanks{Ph.D. Candidate, Kevin T. Crofton Department of Aerospace and Ocean Engineering}, Mayuresh Patil \thanks{Associate Professor, Kevin T. Crofton Department of Aerospace and Ocean Engineering}, and Craig Woolsey \thanks{Professor, Kevin T. Crofton Department of Aerospace and Ocean Engineering}\\ Virginia Polytechnic Institute and State University, Blacksburg, VA, 24061}
\begin{document}
\maketitle
\begin{abstract} 
The goal of this paper is to develop energy-preserving variational integrators for time-dependent mechanical systems with forcing. We first present the Lagrange-d’Alembert principle in the extended Lagrangian mechanics framework and derive the extended forced Euler-Lagrange equations in continuous-time. We then obtain the extended forced discrete Euler-Lagrange equations using the extended discrete mechanics framework and derive adaptive time step variational integrators for time-dependent Lagrangian systems with forcing. We consider two numerical examples to study the numerical performance of energy-preserving variational integrators. First, we consider the example of a nonlinear conservative system to illustrate the advantages of using adaptive time-stepping in variational integrators. We show a trade-off between energy-preserving performance and accurate discrete trajectories when choosing an initial time step. In addition, we demonstrate how the implicit equations become more ill-conditioned as the adaptive time step decreases through a condition number analysis. As a second example, we numerically simulate a damped harmonic oscillator using the adaptive time step variational integrator framework. The adaptive time step increases monotonically for the dissipative system leading to unexpected energy behavior.  

\end{abstract}
%
%
%
%
\section{Introduction} 
In engineering applications, numerical integrators for equations of motion are usually derived by discretizing  differential equations. These traditional integrators do not account for the inherent geometric structure  of the governing continuous-time equations, which results in numerical methods that introduce numerical dissipation and do not preserve invariants of the system.  The field of geometric numerical integration is concerned with numerical methods that preserve the structure of the problem and the corresponding geometric properties of the differential equations. A brief introduction to the field of geometric numerical integration can be found in \cite{Leimkuhler2005SimulatingDynamics} and various techniques for constructing structure-preserving integrators for ordinary differential equations are given in \cite{Hairer2006GeometricEquations} . \par
Ge and Marsden \cite{Zhong1988Lie-PoissonIntegratorsb} showed that a fixed time step numerical integrator cannot preserve the symplectic form, momentum, and energy simultaneously for non-integrable systems. Based on this result, structure-preserving fixed time step mechanical integrators can be divided into two categories, symplectic-momentum, and energy-momentum integrators. Even though symplectic-momentum integrators do not conserve energy exactly, they have been shown to exhibit good long-time energy behavior. On the other hand, Simo and his collaborators \cite{Simo1992TheElastodynamics,Simo1992ExactDynamics} have developed mechanical integrators that conserve energy and momentum but do not preserve the symplectic form. \par
Variational integrators are a class of structure-preserving integrators, that are derived by discretizing the action principle rather than the governing differential equation.  The basic idea behind these methods is to obtain an approximation of the action integral called the discrete action. Stationary points of the discrete action give discrete time trajectories of the mechanical system. These integrators have good long-time energy behavior and they conserve the invariants of the dynamics in the presence of symmetries.  The basics of variational integrators can be reviewed in  \cite{Lew2016AIntegrators}  and more detailed theory can be found in \cite{Marsden2001DiscreteIntegrators}. A more general framework encompassing variational integrators, asynchronous variational integrators, and symplectic-energy-momentum integrators is discussed in \cite{leok2005generalized}. Due to their symplectic nature,  variational integrators are ideal for long-time simulation of conservative or weakly dissipative systems found in astrophysics and molecular dynamics. The discrete trajectories obtained using variational integrators display excellent energy behavior for exponentially long times. \par
The fixed time step variational integrators derived from the discrete variational principle cannot preserve the energy of the system exactly.  To conserve the energy in addition to preserving the symplectic structure and conserving the momentum, time adaptation needs to be used.  These symplectic-energy-momentum integrators were first developed for conservative systems in \cite{Kane1999Symplectic-energy-momentumIntegrators}  by imposing an additional energy preservation equation to compute the time step. The same integrators were derived through a variational approach for a more general case of time-dependent Lagrangian systems in \cite{Marsden2001DiscreteIntegrators} by preserving the discrete energy obtained from the discrete variational principle in the extended phase space. \par
Numerically equivalent methods were derived independently through the Hamiltonian approach by Shibberu in \cite{Shibberu1992Discrete-timeDynamics}. Shibberu has discussed the well-posedness of symplectic-energy-momentum integrators in \cite{Shibberu2006Is} and suggested ways to regularize the governing set of nonlinear discrete equations in \cite{Shibberu2005HowIntegrator}. These symplectic-energy-momentum integrators require solving a coupled nonlinear implicit system of equations at every time step to update the configuration variables and time variable. The time-marching equations for these energy-preserving integrators are ill-conditioned for arbitrarily small time steps and existence of solutions for these discrete trajectories is still an open problem.\par
The purpose of this paper is twofold. First, we develop variational integrators for time-dependent Lagrangian systems with nonconservative forces based on the discretization of the Lagrange-d’Alembert principle in extended phase space. We modify the Lagrange-d’Alembert principle to include time variations in the extended phase space and derive the extended forced Euler-Lagrange equations. We then present a discrete variational principle for time-dependent Lagrangian systems with forcing and derive extended discrete Euler-Lagrange equations. We use the extended discrete mechanics formulation to construct adaptive time step variational integrators for nonautonomous Lagrangian systems with forcing that capture the rate of energy evolution accurately. Second, we consider two numerical examples to understand the numerical properties of the adaptive time step variational integrators and compare the results with fixed time step variational integrators to illustrate the advantages of using adaptive time-stepping in variational integrators.  We also study the effect of initial time step and initial conditions on the numerical performance of the adaptive time step variational integrators. \par
%
%
The remainder of the paper is organized as follows. In Section \ref{s:2}, we review the basics of Lagrangian mechanics and discrete mechanics. We also discretize Hamilton's principle to derive both fixed time step and adaptive time step variational integrators. In Section \ref{s:3}, we modify the Lagrange-d'Alembert principle in the extended phase space and derive extended forced Euler-Lagrange equations. In Section \ref{s:4}, we derive the extended forced discrete Euler-Lagrange equations and obtain adaptive time step variational integrators for time-dependent Lagrangian systems with forcing. In Section \ref{s:5}, we give numerical examples to understand the numerical performance of the adaptive time step variational integrators. Finally, in Section \ref{s:6} we provide concluding remarks and suggest future research directions.
%
%
%
%
\section{Background: Variational Integrators}
\label{s:2} 
In this section, we review the derivation of both fixed and adaptive time step variational integrators. Drawing on the work of Marsden et al. \cite{Kane1999Symplectic-energy-momentumIntegrators, Marsden2001VariationalMechanics, Marsden2001DiscreteIntegrators}, we first derive equations of motion in continuous-time from the variational principle and then derive variational integrators by considering the discretized variational principle in the discrete-time domain.  To this end, we first derive  continuous-time Euler-Lagrange equations of motion from Hamilton's principle. After deriving equations of motion, we use concepts of discrete mechanics developed in \cite{Marsden2001DiscreteIntegrators} to derive discrete Euler-Lagrange equations and then write them in the time-marching form to  obtain both fixed and adaptive time step variational integrators for Lagrangian systems. 
%
%
%
\subsection{Lagrangian Mechanics}
\label{sc:21} 
Hamilton's principle of stationary action is one of the most fundamental results of classical mechanics and is commonly used to derive equations of motion for a variety of systems. The forces and interactions that govern the dynamical evolution of the system are easily determined through Hamilton's principle in a formulaic and elegant manner.  Hamilton's principle \cite{Goldstein1980ClassicalMechanics} states that: The motion of the system between two fixed points from $t_i$ to $t_f$ is such that the action integral has a stationary value for the actual path of the motion. In order to derive the Euler-Lagrange equations via Hamilton's principle, we start by defining the configuration space, tangent space and path space.
 \par
Consider a time-independent Lagrangian system with the  smooth configuration manifold $Q$, state space $TQ$, and a Lagrangian $L: TQ \to \mathbb{R} $. For a given time interval $[t_i,t_f]$, we define the path space to be
\begin{equation}
\mathcal{C}(Q)=\left\{\textbf{q}:[t_i,t_f] \to Q  \ |\  \textbf{q} \ is \  a \ C^2 \  curve\right\}
\end{equation}
where $\textbf{q}$ has the coordinate expression $(q^1,...,q^d)$. The corresponding action $\mathfrak{B}: \mathcal{C}(Q) \to \mathbb{R}$ is obtained by integrating the Lagrangian along the curve $\textbf{q}(t)$, which gives
\begin{equation}
\mathfrak{B}(\textbf{q})=\int^{t_f}_{t_i} L(\textbf{q}(t),\dot{\textbf{q}}(t)) \ dt
\end{equation}
A variation $\delta \textbf{q}(t)$ of $\textbf{q}(t)$ is defined as
\begin{equation}
\delta \textbf{q}(t) = \frac{d}{d\epsilon}\textbf{q}^{\epsilon}(t) \mid _{\epsilon=0}
\end{equation}
where $\textbf{q}^{\epsilon}(t)$ is a smooth one-parameter family of trajectories for some real parameter  $\epsilon$ with $\textbf{q}^0(t)=\textbf{q}(t)$. Hamilton's principle seeks paths $\textbf{q}(t) \in \mathcal{C}(Q)$ which pass through $\textbf{q}(t_i)$ and $\textbf{q}(t_f)$ and also satisfy
\begin{equation}
\delta \mathfrak{B}(\textbf{q})= \delta \ \int^{t_f}_{t_i} L(\textbf{q}(t),\dot{\textbf{q}}(t)) \ dt =0
\end{equation}
We compute variations of the action
\begin{equation}
\delta \mathfrak{B}(\textbf{q})=\int^{t_f}_{t_i}\left[ \frac{\partial L(\textbf{q}(t),\dot{\textbf{q}}(t))}{\partial \textbf{q}} \cdot \delta \textbf{q}  + \frac{\partial L(\textbf{q}(t),\dot{\textbf{q}}(t))}{\partial \dot{\textbf{q}}} \cdot\delta \dot{\textbf{q}} \right] \ dt 
\end{equation}
which after using integration by parts gives
\begin{equation}
\delta \mathfrak{B}(\textbf{q})=\int^{t_f}_{t_i} \left[ \frac{\partial L(\textbf{q}(t),\dot{\textbf{q}}(t))}{\partial \textbf{q}}   - \frac{d}{dt} \left( \frac{\partial L(\textbf{q}(t),\dot{\textbf{q}}(t))}{\partial \dot{\textbf{q}}} \right)   \right]\cdot\delta \textbf{q} \ dt  + \left[ \frac{\partial L(\textbf{q}(t),\dot{\textbf{q}}(t))}{\partial \dot{\textbf{q}}} \cdot \delta \textbf{q} \right] ^{t_f}_{t_i}
\end{equation}
Since the endpoints are fixed, we set the variations of $\textbf{q}(t)$ at the endpoints $\delta \textbf{q}(t_i)$ and $\delta \textbf{q}(t_f)$ equal to zero. Now, requiring that the variations of the action are zero for any admissible $\delta \textbf{q}$ gives the Euler-Lagrange equations for the autonomous Lagrangian system
\begin{equation}
 \frac{\partial L(\textbf{q}(t),\dot{\textbf{q}}(t))}{\partial \textbf{q}}   - \frac{d}{dt} \left( \frac{\partial L(\textbf{q}(t),\dot{\textbf{q}}(t))}{\partial \dot{\textbf{q}}}  \right) =\textbf{0}   
\label{eq:cont}
\end{equation}
which in the coordinate form is 
\begin{equation}
 \frac{\partial L(\textbf{q}(t),\dot{\textbf{q}}(t))}{\partial q^i}   - \frac{d}{dt} \left( \frac{\partial L(\textbf{q}(t),\dot{\textbf{q}}(t))}{\partial \dot{q}^i}  \right) =0    \ \ \ for \ \  i=1,...,d
\end{equation}  \par 
The Lagrange-d'Alembert principle generalizes Hamilton's principle to Lagrangian systems with external forcing. For an autonomous Lagrangian system with forcing (time-independent), we use the Lagrange-d'Alembert principle which seeks paths that satisfy
\begin{equation}
\delta \ \int^{t_f}_{t_i} L(\textbf{q}(t),\dot{\textbf{q}}(t)) \ dt + \int ^{t_f}_{t_i} \textbf{f}_L (\textbf{q}(t),\dot{\textbf{q}}(t)) \cdot \delta \textbf{q} \ dt  =0
\label{eq:fl}
\end{equation}
where the second term accounts for the virtual work done by the forces when the path $\textbf{q}(t)$ is varied by $\delta \textbf{q}(t)$. Using integration by parts gives
\begin{equation}
\int^{t_f}_{t_i} \left[ \frac{\partial L(\textbf{q}(t),\dot{\textbf{q}}(t))}{\partial \textbf{q}}   - \frac{d}{dt} \left( \frac{\partial L(\textbf{q}(t),\dot{\textbf{q}}(t))}{\partial \dot{\textbf{q}}}  \right)  + \textbf{f}_L(\textbf{q}(t),\dot{\textbf{q}}(t))  \right]\cdot\delta \textbf{q} \ dt  + \left[ \frac{\partial L(\textbf{q}(t),\dot{\textbf{q}}(t))}{\partial \dot{\textbf{q}}}\cdot \delta \textbf{q} \right] ^{t_f}_{t_i}=0
\end{equation}
Setting the variations at the endpoints equal to zero, gives the forced Euler-Lagrange equations
\begin{equation}
 \frac{\partial L(\textbf{q}(t),\dot{\textbf{q}}(t))}{\partial \textbf{q}}   - \frac{d}{dt} \left( \frac{\partial L(\textbf{q}(t),\dot{\textbf{q}}(t))}{\partial \dot{\textbf{q}}} \right) + \textbf{f}_L(\textbf{q}(t),\dot{\textbf{q}}(t)) =\textbf{0}
\end{equation}
which can be written in the coordinate form by
\begin{equation}
 \frac{\partial L(\textbf{q}(t),\dot{\textbf{q}}(t))}{\partial q^i}   - \frac{d}{dt} \left( \frac{\partial L(\textbf{q}(t),\dot{\textbf{q}}(t))}{\partial \dot{q}^i} \right) + f^i_L(\textbf{q}(t),\dot{\textbf{q}}(t)) =0  \  \ for \ \  i=1,...,d.
\end{equation}
%
%
%
\subsection{Variational Integrators}
\label{sc:22}  
This subsection presents a systematic approach to derive fixed time step variational integrators for autonomous Lagrangian systems with or without forcing. We use the concepts from discrete Lagrangian mechanics developed in \cite{Marsden2001DiscreteIntegrators} to derive the discrete Euler-Lagrange equations. The discrete trajectories obtained from solving the derived equations are written in a time-marching form to describe the discrete system as an integrator of the continuous system.
\par
Having derived the continuous-time equations of motion from Hamilton's principle, we consider the discrete counterpart to derive fixed time step variational integrators. Consider a discrete Lagrangian system with configuration manifold $Q$ and discrete state space $Q \times Q$. For a fixed time step $h=\frac{t_f-t_i}{N}$, the discrete trajectory is defined by the configuration of the system at the sequence of times $\{\ t_k= t_i + kh|k=0,...,N \}\  \subset \mathbb{R}$. The discrete path space is  
\begin{equation}
\mathcal{C}_d(Q)= \mathcal{C}_d(\{\ t_k \}_{k=0}^N, Q)  = \{\ \textbf{q}_d : \{\ t_k \}_{k=0}^N\ \to Q \}\
\end{equation}
We introduce the discrete Lagrangian function $L_d(\textbf{q}_k,\textbf{q}_{k+1})$, an approximation of the action integral along the curve between $\textbf{q}_k$ and $\textbf{q}_{k+1}$. It approximates the integral of the Lagrangian in the following sense
\begin{equation}
L_d(\textbf{q}_k,\textbf{q}_{k+1}) \approx \int^{t_{k+1}}_{t_k} L(\textbf{q}(t),\dot{\textbf{q}}(t)) \  dt
\end{equation}
The discrete action map $\mathcal{B}_d : \mathcal{C}_d(Q) \to \mathbb{R} $ is defined by 
\begin{equation}
\mathcal{B}_d (\textbf{q}_d) = \sum_{k=0}^{N-1} L_d(\textbf{q}_k,\textbf{q}_{k+1})
\end{equation}
From the discrete Hamilton's principle, we know that  a discrete path $\textbf{q}_d \in \mathcal{C}_d(Q)$ is a solution if the variation of the discrete action is zero for arbitrary variations $\delta \textbf{q}_d \in T_{q_d}\mathcal{C}_d(Q)$ . Taking the variation of the discrete action map gives
\begin{equation}
\delta \mathcal{B}_d(\textbf{q}_d) = \sum_{k=0}^{N-1} [ D_1L_d(\textbf{q}_k,\textbf{q}_{k+1})\cdot \delta \textbf{q}_k  + D_2L_d(\textbf{q}_k,\textbf{q}_{k+1})\cdot\delta \textbf{q}_{k+1} ] 
\end{equation}
where $D_i$ denotes differentiation with respect to the $i^{th}$ argument of the discrete Lagrangian $L_d$. The variation can be rearranged to give
\begin{equation}
\delta \mathcal{B}_d(\textbf{q}_d) = \sum_{k=1}^{N-1} [ D_1L_d(\textbf{q}_k,\textbf{q}_{k+1}) + D_2L_d(\textbf{q}_{k-1},\textbf{q}_{k})]\cdot \delta \textbf{q}_k  + D_1L_d(\textbf{q}_0,\textbf{q}_1)\cdot\delta \textbf{q}_0 + D_2L_d(\textbf{q}_{N-1},\textbf{q}_N)\cdot \delta \textbf{q}_N
\end{equation}
Setting the  variation of the discrete action equal to zero for arbitrary variations gives the discrete Euler-Lagrange equations
\begin{equation}
D_2L_d(\textbf{q}_{k-1},\textbf{q}_{k}) + D_1L_d(\textbf{q}_k,\textbf{q}_{k+1}) =\textbf{0},  \ \ \ \ k=1,..., N-1
\label{eq:ds}
\end{equation}
Given $(\textbf{q}_{k-1},\textbf{q}_k)$, the  above equation can be solved to obtain $(\textbf{q}_k,\textbf{q}_{k+1})$. This defines the discrete Lagrangian map $F_{L_d}: Q \times Q  \to Q \times Q$ and the discrete system \eqref{eq:ds} can be seen as an integrator of the continuous system \eqref{eq:cont}. \par
We define the discrete momentum $\textbf{p}_k=D_2L_d(\textbf{q}_{k-1},\textbf{q}_k)$ for each $k$, which simplifies the second-order discrete time-marching equations \eqref{eq:ds} to first-order form
\begin{equation}
-D_1L_d(\textbf{q}_k,\textbf{q}_{k+1})=\textbf{p}_k
\label{eqn:vi_imp}
\end{equation}
\begin{equation}
\textbf{p}_{k+1}=D_2L_d(\textbf{q}_{k},\textbf{q}_{k+1})
\label{eqn:vi_exp}
\end{equation}
Given $(\textbf{q}_k,\textbf{p}_k)$, \eqref{eqn:vi_imp} is solved implicitly to obtain $\textbf{q}_{k+1}$, which is then used to obtain $\textbf{p}_{k+1}$ from \eqref{eqn:vi_exp}. \par
For cases with forcing, we define two discrete forces $\textbf{f}_d^{\pm} :  Q \times Q \to T^*Q$ which approximate the continuous-time force integral in the following sense
\begin{equation}
 \textbf{f}_d^+(\textbf{q}_k,\textbf{q}_{k+1})\cdot\delta \textbf{q}_{k+1} +  \textbf{f}_d^-(\textbf{q}_k,\textbf{q}_{k+1})\cdot\delta \textbf{q}_{k}  \approx \int ^{t_{k+1}}_{t_k} \textbf{f}_L (\textbf{q}(t), \dot{\textbf{q}}(t)) .\cdot\delta \textbf{q} \ dt
\end{equation}
The discrete Lagrange-d'Alembert principle seeks curves $\{\ \textbf{q}_k  \  \}_{k=0}^N\ $ that satisfy 
\begin{equation}
\delta \sum_{k=0}^{N-1} L_d(\textbf{q}_k,\textbf{q}_{k+1}) + \sum_{k=0}^{N-1}[\textbf{f}_d^+(\textbf{q}_k,\textbf{q}_{k+1})\cdot\delta \textbf{q}_{k+1} +  \textbf{f}_d^-(\textbf{q}_k,\textbf{q}_{k+1})\cdot\delta \textbf{q}_{k}] = 0
\end{equation}
which gives the forced discrete Euler-Lagrange equations 
\begin{equation}
D_2L_d(\textbf{q}_{k-1},\textbf{q}_{k}) + D_1L_d(\textbf{q}_k,\textbf{q}_{k+1}) + \textbf{f}_d^+(\textbf{q}_{k-1},\textbf{q}_k) + \textbf{f}_d^-(\textbf{q}_{k},\textbf{q}_{k+1}) = \textbf{0}   \ \ \ \ k=1,..., N-1
\end{equation}
These equations can be implemented as a variational integrator for forced Lagrangian systems as follows
\begin{equation}
-D_1L_d(\textbf{q}_k,\textbf{q}_{k+1}) - \textbf{f}_d^-(\textbf{q}_k,\textbf{q}_{k+1}) =\textbf{p}_k
\label{eqn:vif_imp}
\end{equation}
\begin{equation}
\textbf{p}_{k+1}=D_2L_d(\textbf{q}_{k},\textbf{q}_{k+1}) + \textbf{f}_d^+(\textbf{q}_k,\textbf{q}_{k+1})
\label{eqn:vif_exp}
\end{equation}
\par 
The fixed time step variational integrators \eqref{eqn:vi_imp}-\eqref{eqn:vi_exp} are symplectic and momentum-preserving. These fixed time step algorithms have bounded energy error; the magnitude of the energy error depends on where the trajectory is in the phase space. To preserve the energy of the Lagrangian system exactly, time adaptation must be used. To this end, we use the extended Lagrangian mechanics framework where, in addition to dynamic configuration variables, time is also treated as a dynamic variable.
%
%
%
\subsection{Extended Lagrangian Mechanics}
\label{sc:23}

Consider a time-dependent Lagrangian system with configuration manifold $Q$ and time space $\mathbb{R}$. In the extended Lagrangian mechanics framework \cite{Marsden2001DiscreteIntegrators}, we treat time as a dynamic variable and define the extended configuration manifold $\bar{Q}=\mathbb{R} \times Q$; the corresponding state space $T\bar{Q}$ is $\mathbb{R} \times TQ$. The extended Lagrangian is $L: \mathbb{R} \times TQ \to \mathbb{R}$. \par
In the extended Lagrangian mechanics framework, $t$ and $\textbf{q}$ are both parametrized by an independent variable $a$. The two components of a trajectory $c$ are $c(a)=(c_t(a),c_{\textbf{q}}(a))$. The extended path space is 
\begin{equation}
\bar{\mathcal{C}}=\left\{c:[a_0,a_f] \to \bar{Q} |\  c \ is \  a \ C^2 \  curve \ and \ c'_t(a) >0 \right\}
\end{equation}
For a given path $c(a)$, the initial time is $t_0=c_t(a_0)$ and the final time is $t_f=c_t(a_f)$. The extended action $\bar{\mathfrak{B}} : \bar{\mathcal{C}} \to \mathbb{R}$ is
\begin{equation}
\bar{\mathfrak{B}} = \int ^{t_f} _{t_0} L(t, \textbf{q}(t), \dot{\textbf{q}}(t)) dt
\end{equation}
Since time is a dynamic variable in this framework, we substitute $(t,\textbf{q}(t),\dot{\textbf{q}}(t))=\left(c_t(a),c_{\textbf{q}}(a),\frac{c'_{\textbf{q}}(a)}{c'_t(a)} \right)$ in the above equation to get 
\begin{equation}
\bar{\mathfrak{B}} = \int ^{a_f} _{a_0} L\left(c_t(a),c_q(a),\frac{c'_q(a)}{c'_t(a)} \right)c'_t(a) \  da
\end{equation}
We compute variations of the action
\begin{equation}
\delta \bar{\mathfrak{B}} = \int^{a_f}_{a_0} \left[ \frac{\partial L}{\partial t} \delta c_t + \frac{\partial L}{\partial \textbf{q}}\cdot \delta c_{\textbf{q}}  +  \frac{\partial L}{\partial \dot{\textbf{q}}}\cdot \left( \frac{\delta c'_{\textbf{q}}(a)}{c'_t}  -\frac{c'_\textbf{q} \delta c'_t(a) }{(c'_t)^2} \right) \right] c'_t(a)da + \int ^{a_f} _{a_0} L \delta c'_t(a) \  da
\end{equation}
Using integration by parts and setting the variations at the end points to zero gives
\begin{equation}
\delta \bar{\mathfrak{B}} = \int^{a_f}_{a_0} \left[ \frac{\partial L}{\partial \textbf{q}} c'_t - \frac{d}{da} \frac{\partial L}{\partial \dot{\textbf{q}}} \right] \cdot \delta c_\textbf{q}(a)da + \int ^{a_f} _{a_0}  \left[ \frac{\partial L}{\partial t}c'_t + \frac{d}{da}\left( \frac{\partial L}{ \partial \dot{\textbf{q}}}. \frac{c'_\textbf{q}}{c'_t} - L \right) \right] \delta c_t(a) \  da
\end{equation}
Using $dt=c'_t(a) da$ in the above expression gives two equations of motion. The first is the Euler-Lagrange equation of motion 
\begin{equation}
\frac{\partial L}{\partial \textbf{q}}   - \frac{d}{dt} \left( \frac{\partial L}{\partial \dot{\textbf{q}}}  \right) =\textbf{0}  
\label{eq:contp}
\end{equation}
which is the same as the equation obtained using the classical Lagrangian mechanics framework in Section \ref{sc:21}. The second equation is  
\begin{equation}
\frac{\partial L}{ \partial t} + \frac{d}{dt} \left( \frac{\partial L}{ \partial \dot{\textbf{q}}} . \dot{\textbf{q}} - L \right) =0
\label{eq:conte}
\end{equation}
which describes how the energy of the system evolves with time. 
\subsection{Energy-preserving Variational Integrators}
\label{sc:24}  
For the extended Lagrangian mechanics, we define the extended discrete state space $\bar{Q} \times \bar{Q}$. The extended discrete path space is 
\begin{equation}
\bar{\mathcal{C}}_d = \{\ c: \{\ 0,...,N \} \to \bar{Q} | \ \  c_t(k+1) > c_t(k) \ \ for \  all \  k  \}\
\end{equation}
The extended discrete action map $\bar{\mathcal{B}}_d : \bar{\mathcal{C}}_d \to \mathbb{R}$ is
\begin{equation}
\bar{\mathcal{B}}_d = \sum _{k=0}^{N-1} L_d(t_k,\textbf{q}_k,t_{k+1},\textbf{q}_{k+1})
\end{equation}
where $L_d : \bar{Q} \times \bar{Q} \to \mathbb{R}$ is the extended discrete Lagrangian function which approximates the action integral between two successive configurations. Taking variations of the extended discrete action map gives
\begin{multline}
      \delta \bar{B}_d = \sum_{k=1}^{N-1} \left[ D_4L_d(t_{k-1},\textbf{q}_{k-1},t_{k},\textbf{q}_{k}) + D_2L_d(t_k,\textbf{q}_k,t_{k+1},\textbf{q}_{k+1}) \right] \cdot \delta \textbf{q}_k  \\  +\sum_{k=1}^{N-1} \left[ D_3L_d(t_{k-1},\textbf{q}_{k-1},t_{k},\textbf{q}_{k}) + D_1L_d(t_k,\textbf{q}_k,t_{k+1},\textbf{q}_{k+1}) \right]  \delta t_k  = 0
    \end{multline}
Applying Hamilton's principle of least action and setting variations at end points to zero gives the extended discrete Euler-Lagrange equations 
\begin{equation}
 D_4L_d(t_{k-1},\textbf{q}_{k-1},t_{k},\textbf{q}_{k}) + D_2L_d(t_k,\textbf{q}_k,t_{k+1},\textbf{q}_{k+1}) = \textbf{0}
\label{eq:vi_pvv}
\end{equation}

\begin{equation}
D_3L_d(t_{k-1},\textbf{q}_{k-1},t_{k},\textbf{q}_{k}) + D_1L_d(t_k,\textbf{q}_k,t_{k+1},\textbf{q}_{k+1}) = 0
\label{eq:vi_ev}
\end{equation}
Given $(t_{k-1},\textbf{q}_{k-1},t_{k},\textbf{q}_{k})$, the extended discrete Euler-Lagrange equations can be solved to obtain $\textbf{q}_{k+1}$ and $t_{k+1}$. This  extended discrete Lagrangian system can be seen as a numerical integrator of the continuous-time nonautonomous Lagrangrian system with adaptive time steps.
\par
In the extended discrete mechanics framework, we define the discrete momentum $\textbf{p}_k$ by
\begin{equation}
\textbf{p}_{k}= D_4L_d(t_{k-1},\textbf{q}_{k-1},t_{k},\textbf{q}_{k})
\end{equation}
We also introduce the discrete energy 
\begin{equation}
E_{k}=D_3L_d(t_{k-1},\textbf{q}_{k-1},t_{k},\textbf{q}_{k})
\end{equation}
Using the discrete momentum and discrete energy definitions, we can re-write the extended discrete Euler-Lagrange equations \eqref{eq:vi_pvv} and \eqref{eq:vi_ev} in the following form
\begin{equation}
-D_2L_d(t_k,\textbf{q}_k,t_{k+1},\textbf{q}_{k+1}) = \textbf{p}_k
\label{eq:pimp}
\end{equation}
\begin{equation}
D_1L_d(t_k,\textbf{q}_k,t_{k+1},\textbf{q}_{k+1}) = E_k
\label{eq:eimp}
\end{equation}
\begin{equation}
\textbf{p}_{k+1}= D_4L_d(t_k,\textbf{q}_k,t_{k+1},\textbf{q}_{k+1})
\label{eq:pexp}
\end{equation}
\begin{equation}
E_{k+1}= -D_3L_d(t_k,\textbf{q}_k,t_{k+1},\textbf{q}_{k+1})
\label{eq:eexp}
\end{equation}
Given $(t_k,\textbf{q}_k,\textbf{p}_k,E_k)$, the coupled nonlinear equations \eqref{eq:pimp} and \eqref{eq:eimp} are solved implicitly to obtain $\textbf{q}_{k+1}$ and $t_{k+1}$. The configuration $\textbf{q}_{k+1}$ and time $t_{k+1}$ are then used in \eqref{eq:pexp} and \eqref{eq:eexp} to obtain $(\textbf{p}_{k+1},E_{k+1})$ explicitly. The extended discrete Euler-Lagrange equations were first written in the  time-marching form in \cite{Marsden2001DiscreteIntegrators} and are also known as symplectic-energy-momentum integrators. 
%
%
%
%
%
\section{Modified Lagrange-d'Alembert Principle}
\label{s:3} 
The extended discrete Euler-Lagrange equations derived in Section \ref{sc:24} can be used as energy-preserving variational integrators for Lagrangian systems. In order to extend this energy-preserving variational integrator framework to Lagrangian systems with external forcing, we need to discretize the Lagrange-d'Alembert principle in the extended Lagrangian mechanics framework. We first present  the Lagrange- d'Alembert principle in extended phase space for time-dependent Lagrangian systems with forcing by considering the variations with respect to time $t$. Using the extended Lagrangian mechanics framework, we derive the extended Euler-Lagrange equations for time-dependent Lagrangian systems with forcing.
\par 
As mentioned in  Section \ref{sc:21}, the Lagrange-  d'Alembert principle  \eqref{eq:fl}   modifies Hamilton's principle of stationary action by considering the virtual work done by the forces for a variation $\delta \textbf{q}$ in the configuration variable $\textbf{q}$. Since the standard Lagrangian mechanics framework treats time only as an independent continuous parameter, it does not account for time variations in the Lagrange-d'Alembert principle. Thus, we need to modify the Lagrange-d'Alembert principle in the extended Lagrangian mechanics framework to account for time variations.\par
We modify the Lagrange-d'Alembert principle by adding an additional term in the variational principle that accounts for virtual work done by the external force $\textbf{f}_L$ due to variations in the time variable
\begin{equation}
\delta \ \int^{t_f} _{t_0} L(t,\textbf{q}(t),\dot{\textbf{q}}(t)) \ dt + \int ^{t_f} _{t_0} \textbf{f}_L (t,\textbf{q}(t), \dot{\textbf{q}}(t)) \cdot \delta \textbf{q} \ dt  - \int ^{t_f} _{t_0} \textbf{f}_L (t,\textbf{q}(t), \dot{\textbf{q}}(t)) \cdot \left(\dot{\textbf{q}} \delta t \right) \ dt =0
\end{equation}
Using the extended Lagrangian mechanics framework discussed in Section \ref{sc:23} to derive the equations of motion, we first re-write the modified Lagrange-d'Alembert principle 
\begin{equation}
\delta \int ^{a_f} _{a_0} L\left(c_t(a),c_{\textbf{q}}(a),\frac{c'_{\textbf{q}}(a)}{c'_t(a)} \right)c'_t(a) \  da + \int ^{a_f} _{a_0} f_L\left(c_t(a),c_{\textbf{q}}(a),\frac{c'_{\textbf{q}}(a)}{c'_t(a)} \right)\cdot(\delta c_{\textbf{q}} - \dot{\textbf{q}} \delta c_t ) \  c't(a) da =0
\end{equation}
Taking variations of the discrete action with respect to both configuration $\textbf{q}$ and time $t$ gives
\begin{multline}
 \int^{a_f}_{a_0} \left[ \frac{\partial L}{\partial t} \delta c_t + \frac{\partial L}{\partial \textbf{q}}\cdot \delta c_{\textbf{q}}  +  \frac{\partial L}{\partial \dot{\textbf{q}}}\cdot \left( \frac{\delta c'_{\textbf{q}}(a)}{c'_t}  -\frac{c'_\textbf{q} \delta c'_t(a) }{(c'_t)^2} \right) \right] c'_t(a)da + \int ^{a_f} _{a_0} L \delta c'_t(a) \  da  \\  +  \int ^{a_f} _{a_0} \left( f_L\cdot\delta c_{\textbf{q}}\right) \  c't(a) da - \int ^{a_f} _{a_0}  \left( f_L\cdot\dot{\textbf{q}} \right) \delta c_t  \  c't(a) da  =0
\end{multline}
Using integration by parts and setting variations at the end points to zero gives
\begin{equation}
\int^{a_f}_{a_0} \left[ \frac{\partial L}{\partial \textbf{q}} c'_t - \frac{d}{da} \frac{\partial L}{\partial \dot{\textbf{q}}} + \textbf{f}_L \right] \cdot\delta c_{\textbf{q}}(a)da + \int ^{a_f} _{a_0}  \left[ \frac{\partial L}{\partial t}c'_t + \frac{d}{da}\left( \frac{\partial L}{ \partial \dot{\textbf{q}}} .\frac{c'_{\textbf{q}}}{c'_t} - L \right) - \textbf{f}_L\cdot\dot{\textbf{q}} \right] \delta c_t(a) \  da = 0
\end{equation}
Using $dt=c'_t(a) da$ in the above expression gives two equations of motion for the forced time-dependent Lagrangian system. The first equation is the well-known forced Euler-Lagrange equation for a time-dependent system
\begin{equation}
 \frac{\partial L}{\partial \textbf{q}}   - \frac{d}{dt} \left( \frac{\partial L}{\partial \dot{\textbf{q}}}  \right) + \textbf{f}_L =\textbf{0} 
 \label{eq:fel}
\end{equation}
whereas the second equation is the energy evolution equation 
\begin{equation}
\frac{\partial L}{ \partial t} + \frac{d}{dt} \left( \frac{\partial L}{ \partial \dot{\textbf{q}}} \dot{\textbf{q}} - L \right) - \textbf{f}_L\cdot\dot{\textbf{q}} =0
\label{eq:fev}
\end{equation}
Thus, for the forced case, the energy evolution equation describes how the energy of the Lagrangian system depends on the input power by the external force $\textbf{f}_L$. If we consider an associated curve $\textbf{q}(t)$ satisfying the forced Euler-Lagrange equations and compute the energy evolution equation we get
\begin{equation}
\frac{\partial L}{ \partial t} + \frac{d}{dt} \left( \frac{\partial L}{ \partial \dot{\textbf{q}}} . \dot{\textbf{q}} - L \right) - \textbf{f}_L.\dot{\textbf{q}} = \frac{\partial L}{ \partial t} +    \frac{d}{dt} \left( \frac{\partial L}{ \partial \dot{\textbf{q}}} \right). \dot{\textbf{q}} + \frac{\partial L}{ \partial \dot{\textbf{q}}}. \ddot{\textbf{q}} - \frac{dL}{dt} - \textbf{f}_L.\dot{\textbf{q}}
\end{equation}
which after substituting $\frac{d}{dt} \left ( \frac{\partial L}{ \partial \dot{\textbf{q}}} \right) = \frac{\partial L}{\partial \textbf{q}}  + \textbf{f}_L $ simplifies to 
\begin{equation}
\frac{\partial L}{ \partial t} + \frac{d}{dt} \left( \frac{\partial L}{ \partial \dot{\textbf{q}}} . \dot{\textbf{q}} - L \right) - \textbf{f}_L.\dot{\textbf{q}}= \left (  \frac{\partial L}{ \partial t} +  
\left( \frac{\partial L}{\partial \textbf{q}} + \textbf{f}_L \right) . \dot{\textbf{q}} + \frac{\partial L}{ \partial \dot{\textbf{q}}}. \ddot{\textbf{q}}  \right) - \frac{dL}{dt} -  \textbf{f}_L.\dot{\textbf{q}} =0
\end{equation}
which shows that \eqref{eq:fel} implies \eqref{eq:fev}. Thus,  for continuous-time forced Lagrangian systems, the additional energy evolution equation obtained by taking the variation with respect to time does not provide any new information concerning the forced Euler-Lagrange equations. \\ \\
\textbf{Remark 1.} It should be noted that both \eqref{eq:fel} and \eqref{eq:fev} depend only on the associated curve $\textbf{q}(t)$ and the time component $c_t(a)$ of the extended path cannot be determined from the governing equations. Thus, the "velocity" of time, i.e. $c'_t(a)$, is indeterminate in the continuous-time  formulation of time-dependent Lagrangian systems with forcing. 
%
%
%
%
\section{Energy-preserving Variational Integrators}
\label{s:4} 
In this section, we derive the extended discrete Euler-Lagrange equations for time-dependent Lagrangian systems with forcing by discretizing the modified Lagrange-d'Alembert principle given in Section \ref{s:3}. The key difference from the extended discrete Euler-Lagrange equations derived in Section \ref{sc:24} is that we will have additional discrete terms accounting for the virtual work done by the external forcing. 
\par
The modified Lagrange-d'Alembert principle presented in Section \ref{s:3} has two continuous-time force integrals in the variational principle. In order to derive the extended forced discrete Euler-Lagrange equations, we define two discrete force terms  $f_d^{\pm} :  \bar{Q} \times \bar{Q} \to T^*\bar{Q}$ which approximate the virtual work done due to variations in $\textbf{q}$ in the following sense
\begin{equation}
 f_d^+(t_k,\textbf{q}_k,t_{k+1},\textbf{q}_{k+1})\cdot\delta \textbf{q}_{k+1} +  f_d^-(t_k,\textbf{q}_k,t_{k+1},\textbf{q}_{k+1})\cdot\delta \textbf{q}_{k}  \approx \int ^{t_{k+1}}_{t_k} \textbf{f}_L (t,\textbf{q}(t), \dot{\textbf{q}}(t)) \cdot \delta \textbf{q} \ dt
\end{equation}
 We also define two discrete power terms  $g_d^{\pm} :  \bar{Q} \times \bar{Q} \to \mathbb{R}$ which approximate the virtual work done due to time variations in the following sense
 \begin{equation}
 g_d^+(t_k,\textbf{q}_k,t_{k+1},\textbf{q}_{k+1})\delta t_{k+1} +  g_d^-(t_k,\textbf{q}_k,t_{k+1},\textbf{q}_{k+1})\delta t_{k}  \approx \int ^{t_{k+1}}_{t_k} -\textbf{f}_L (t,\textbf{q}(t), \dot{\textbf{q}}(t)) \cdot \left(\dot{\textbf{q}} \delta t \right) \ dt
\end{equation}
For the time-dependent Lagrangian system with forcing, we seek discrete-time paths which satisfy 
\begin{multline}
      \delta \sum _{k=0}^{N-1} L_d(t_k,\textbf{q}_k,t_{k+1},\textbf{q}_{k+1}) +  \sum_{k=0}^{N-1}[\textbf{f}_d^+(t_k,\textbf{q}_k,t_{k+1},\textbf{q}_{k+1})\cdot\delta \textbf{q}_{k+1} +  \textbf{f}_d^-(t_k,\textbf{q}_k,t_{k+1},\textbf{q}_{k+1})\cdot\delta \textbf{q}_{k}] \\  + \sum_{k=0}^{N-1}[g_d^+(t_k,\textbf{q}_k,t_{k+1},\textbf{q}_{k+1})\delta t_{k+1} +  g_d^-(t_k,\textbf{q}_k,t_{k+1},\textbf{q}_{k+1})\delta t_{k}]  = 0
      \label{eq:ftnew}
    \end{multline}
Setting all the variations at the endpoints equal to zero in \eqref{eq:ftnew} gives the  extended forced  discrete Euler-Lagrange equations
    %
    %
\begin{equation}
D_4L_d(t_{k-1},\textbf{q}_{k-1},t_{k},\textbf{q}_{k}) + \  D_2L_d(t_k,\textbf{q}_k,t_{k+1},\textbf{q}_{k+1}) +  \ \textbf{f}_d^+(t_{k-1},\textbf{q}_{k-1},t_{k},\textbf{q}_{k}) +  \  \textbf{f}_d^-(t_k,\textbf{q}_k,t_{k+1},\textbf{q}_{k+1}) = \textbf{0}
\label{eq:vi_mv_f}
\end{equation}
\begin{equation}
D_3L_d(t_{k-1},\textbf{q}_{k-1},t_{k},\textbf{q}_{k}) + \  D_1L_d(t_k,\textbf{q}_k,t_{k+1},\textbf{q}_{k+1}) + \  g_d^+(t_{k-1},\textbf{q}_{k-1},t_{k},\textbf{q}_{k}) + \  g_d^-(t_k,\textbf{q}_k,t_{k+1},\textbf{q}_{k+1})  = 0
\label{eq:vi_ev_f}
\end{equation}
We modify the definitions of the discrete momentum and energy to account for the effect of forcing
\begin{equation}
\textbf{p}_k=  D_4L_d(t_{k-1},\textbf{q}_{k-1},t_{k},\textbf{q}_{k}) +  \textbf{f}_d^+(t_{k-1},\textbf{q}_{k-1},t_{k},\textbf{q}_{k})
\label{eq: pkmdef}
\end{equation}
\begin{equation}
E_k= - D_3L_d(t_{k-1},\textbf{q}_{k-1},t_{k},\textbf{q}_{k}) -  g_d^+(t_{k-1},\textbf{q}_{k-1},t_{k},\textbf{q}_{k})
\label{eq: ekmdef}
\end{equation}
Using the modified discrete momentum and energy definitions \eqref{eq: pkmdef} and \eqref{eq: ekmdef}, the extended forced discrete Euler-Lagrange equations can be re-written in the following form
\begin{equation}
-D_2L_d(t_k,\textbf{q}_k,t_{k+1},\textbf{q}_{k+1}) - \textbf{f}_d^-(t_k,\textbf{q}_k,t_{k+1},\textbf{q}_{k+1}) = \textbf{p}_k
\label{eq:pimpf}
\end{equation}
\begin{equation}
D_1L_d(t_k,\textbf{q}_k,t_{k+1},\textbf{q}_{k+1}) + g_d^-(t_k,\textbf{q}_k,t_{k+1},\textbf{q}_{k+1}) = E_k
\label{eq:eimpf}
\end{equation}
\begin{equation}
\textbf{p}_{k+1}= D_4L_d(t_k,\textbf{q}_k,t_{k+1},\textbf{q}_{k+1}) + \textbf{f}_d^+(t_k,q_k,t_{k+1},q_{k+1})
\label{eq:pexpf}
\end{equation}
\begin{equation}
E_{k+1}= -D_3L_d(t_k,\textbf{q}_k,t_{k+1},\textbf{q}_{k+1}) - g_d^-(t_k,\textbf{q}_k,t_{k+1},\textbf{q}_{k+1})
\label{eq:eexpf}
\end{equation}
Given a time-dependent Lagrangian system with external forcing, the  extended discrete Lagrangian system obtained by solving \eqref{eq:pimpf}-\eqref{eq:eexpf} can be used as an adaptive time step variational integrator for the continuous-time system.
\\ \\
\textbf{Remark 2.} The modified discrete energy \eqref{eq: ekmdef} has a contribution from the external forcing which accounts for the virtual work done during the adaptive time step. Thus, the discrete trajectory obtained by solving the extended discrete Euler-Lagrange equations preserves a discrete quantity which is not the discrete analogue of the total energy of the Lagrangian system. This detail becomes important when we simulate a dissipative Lagrangian system with an adaptive time step variational integrator in Section \ref{sc:52}
%
%
%
%
\section{Numerical Examples}
\label{s:5}
In this section, we implement the extended discrete Euler-Lagrange equations as numerical integrators of continuous dynamical systems. We  first consider a  nonlinear conservative dynamical system studied in \cite{Kane1999Symplectic-energy-momentumIntegrators} and compare the fixed time step variational integrator results with the corresponding results for  the adaptive time step variational integrator. We then study the damped harmonic oscillator, a time-independent dynamical system, using the extended forced discrete Euler-Lagrange equations in order to investigate the numerical properties of the adaptive time step variational integrators for forced systems.
\subsection{Conservative Example}
\label{sc:51}
We consider a particle in a double-well potential. The Lagrangian for this conservative one degree of freedom dynamical system is 
\begin{equation}
L(q,\dot{q}) = \frac{1}{2}m\dot{q}^2 - V(q)
\end{equation}
where
\begin{equation}
V(q)=\frac{1}{2} \left( q^4 -q^2\right)
\end{equation}
\subsubsection{Fixed time step algorithm}
For the fixed time step case, we choose a constant time step $h$. The discrete Lagrangian is obtained using the midpoint rule
\begin{equation}
L_d(q_k,q_{k+1})=hL\left (\frac{q_k + q_{k+1}}{2},\frac{q_{k+1}-q_k}{h}\right)
\end{equation}
Using \eqref{eqn:vi_exp}, the discrete momentum $p_{k+1}$ is given by \begin{equation}
p_{k+1}=m\left(\frac{q_{k+1}-q_k}{h} \right) + h \left( \left( \frac{q_{k+1}+q_k}{4} \right) -  \left( \frac{q_{k+1}+q_k}{2} \right)^3\right) 
\label{eq:vi_exp_exa}
\end{equation}
For given $(q_k,p_k)$ at the $k^{th}$ time step, using \eqref{eqn:vi_imp} gives the following implicit equation
\begin{equation}
m\left(\frac{q_{k+1}-q_k}{h} \right) - h \left( \left( \frac{q_{k+1}+q_k}{4} \right) -  \left( \frac{q_{k+1}+q_k}{2} \right)^3\right) =p_k
\end{equation}
\begin{figure}[h]
    \centering
    \begin{subfigure}[b]{0.45\textwidth}
        \includegraphics[width=\textwidth]{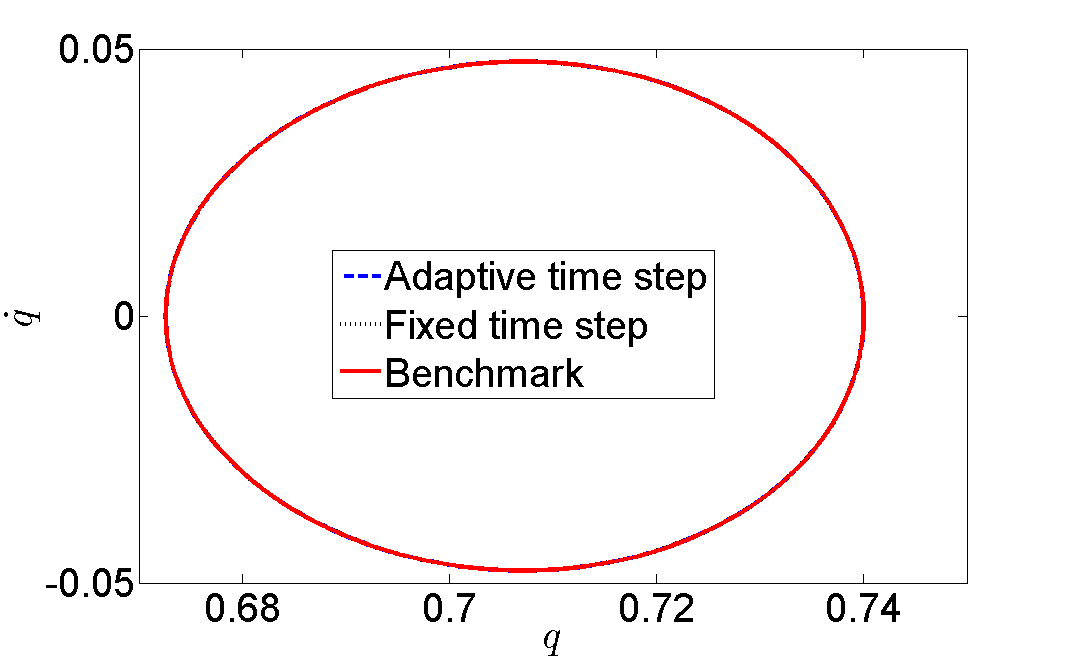}
        \caption{$q(0)=0.74, \ \dot{q}(0)=0$}
        \label{fig:gull0}
    \end{subfigure}
    \begin{subfigure}[b]{0.45\textwidth}
        \includegraphics[width=\textwidth]{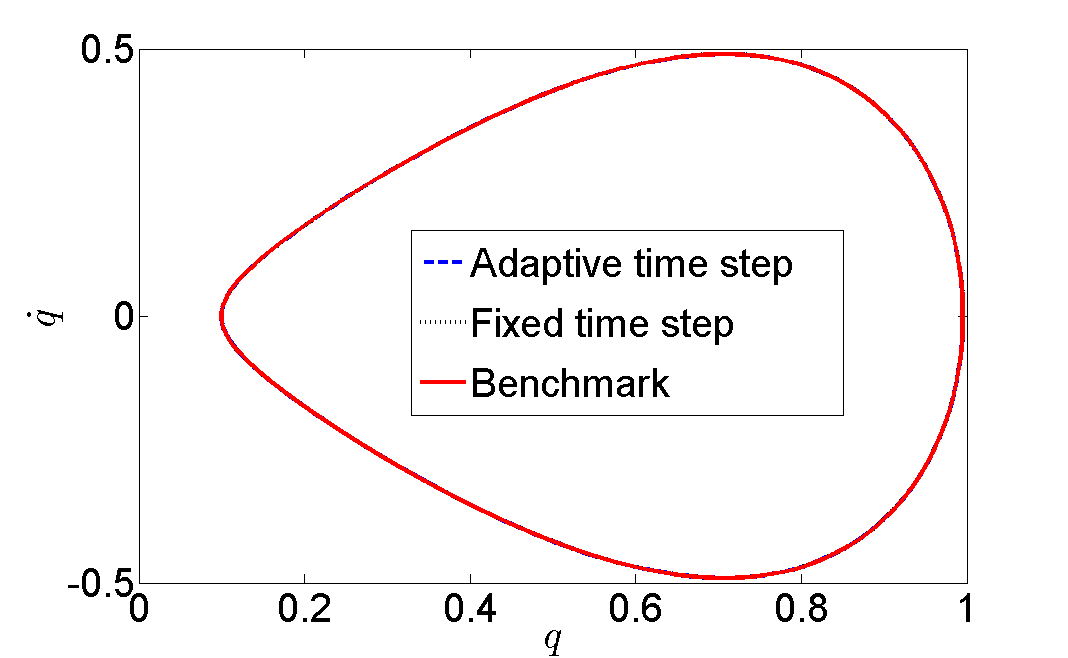}
        \caption{$q(0)=0.995, \ \dot{q}(0)=0$}
        \label{fig:tiger0}
    \end{subfigure}
    \caption{Two initial conditions are studied for the energy error analysis for the particle in double-well potential. An initial time step of $h_0=0.01$ is used for both cases and phase space trajectories for both fixed time step and adaptive time step algorithms are compared to the benchmark trajectory. The trajectories in each figure are indistinguishable.}\label{fig:animals0}
\end{figure} 
\subsubsection{Adaptive time step algorithm}
For the adaptive time step case, we have discrete time $t_k$  as an additional discrete variable. We use the midpoint rule to obtain the discrete Lagrangian $L_d$
\begin{equation}
L_d(t_k,q_k,t_{k+1},q_{k+1})=(t_{k+1}-t_k) \ \left[ \frac{1}{2}m \left( \frac{q_{k+1}-q_k}{t_{k+1}-t_k} \right)^2 - \frac{1}{2} \left( \left( \frac{q_{k+1}+q_k}{2} \right)^4 -\left( \frac{q_{k+1}+q_k}{2} \right)^2\right) \right ] 
\label{eq:adap}
\end{equation}
The discrete momentum $p_k$ and discrete energy $E_k$ are obtained by substituting the $L_d$ expression in \eqref{eq:pexp} and \eqref{eq:eexp}
\begin{equation}
p_{k+1} =  m\left(\frac{q_{k+1}-q_k}{t_{k+1}-t_k} \right) + (t_{k+1}-t_k) \left( \left( \frac{q_{k+1}+q_k}{4} \right) -  \left( \frac{q_{k+1}+q_k}{2} \right)^3\right)
\end{equation}
\begin{equation}
E_{k+1}= \frac{1}{2}m \left( \frac{q_{k+1}-q_k}{t_{k+1}-t_k} \right)^2 +\  \frac{1}{2} \left( \left( \frac{q_{k+1}+q_k}{2} \right)^4 -\left( \frac{q_{k+1}+q_k}{2} \right)^2\right)
\label{eq:disce}
\end{equation}
The implicit time-marching equations for the adaptive time step algorithm are 
\begin{equation}
m\left(\frac{q_{k+1}-q_k}{t_{k+1}-t_k} \right) - (t_{k+1}-t_k) \left( \left( \frac{q_{k+1}+q_k}{4} \right) -  \left( \frac{q_{k+1}+q_k}{2} \right)^3\right)=p_k
\end{equation}
\begin{equation}
\frac{1}{2}m \left( \frac{q_{k+1}-q_k}{t_{k+1}-t_k} \right)^2 +\  \frac{1}{2} \left( \left( \frac{q_{k+1}+q_k}{2} \right)^4 -\left( \frac{q_{k+1}+q_k}{2} \right)^2\right)=E_k
\end{equation}
Since the dynamical system being considered here is time-independent, we rewrite the time-marching equations in terms of $h_k=(t_{k+1}-t_k)$ and $v_k= \left( \frac{q_{k+1}-q_k}{t_{k+1}-t_k} \right)$
\begin{equation}
F(q_k,p_k,h_k,v_k)=mv_k -h_k \left( \left( \frac{v_kh_k+ 2q_k}{4} \right) -  \left( q_k + \frac{v_kh_k}{2} \right)^3\right)-p_k=0
\label{eq:tmpc}
\end{equation}
\begin{equation}
G(q_k,E_k,h_k,v_k)= \frac{1}{2}mv_k^2 + \frac{1}{2} \left( \left( q_k + \frac{v_kh_k}{2} \right)^4 - \left( q_k + \frac{v_kh_k}{2} \right)^2\right) - E_k =0 
\label{eq:tmec}
\end{equation}
These time-marching equations are solved using Newton's iterative method with the restriction $h_k>0$ to obtain discrete trajectories in the extended space. This extended discrete system can be used as a variational integrator for the continuous-time dynamical system.
\begin{figure}[h]
    \centering
    \begin{subfigure}[b]{0.45\textwidth}
        \includegraphics[width=\textwidth]{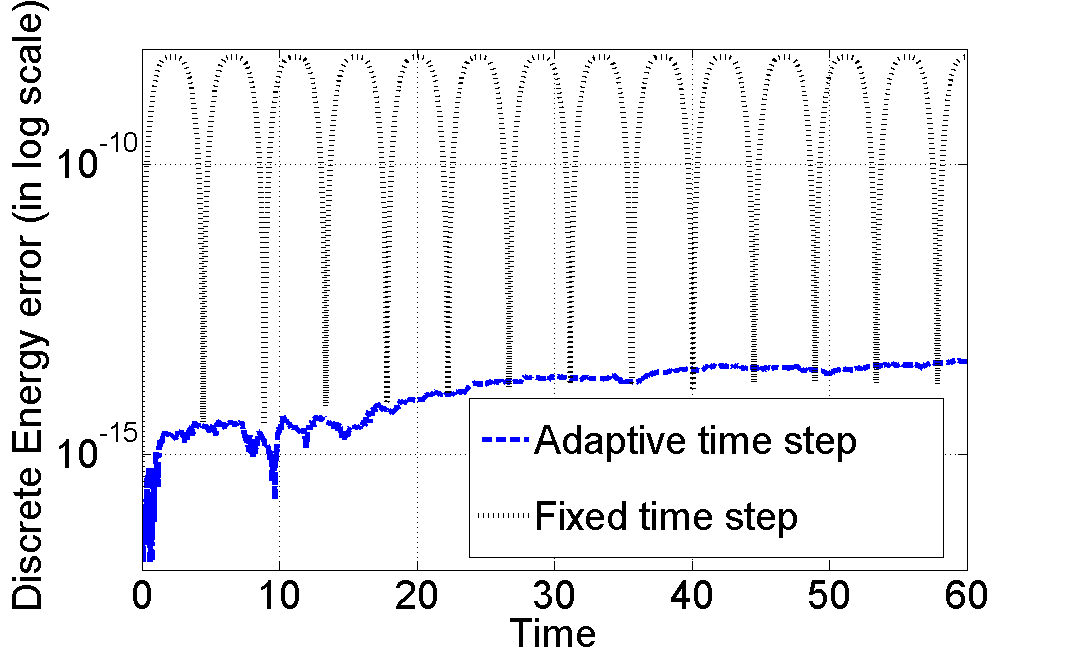}
        \caption{$q(0)=0.74, \ \dot{q}(0)=0$}
        \label{fig:gull1}
    \end{subfigure}
    \begin{subfigure}[b]{0.45\textwidth}
        \includegraphics[width=\textwidth]{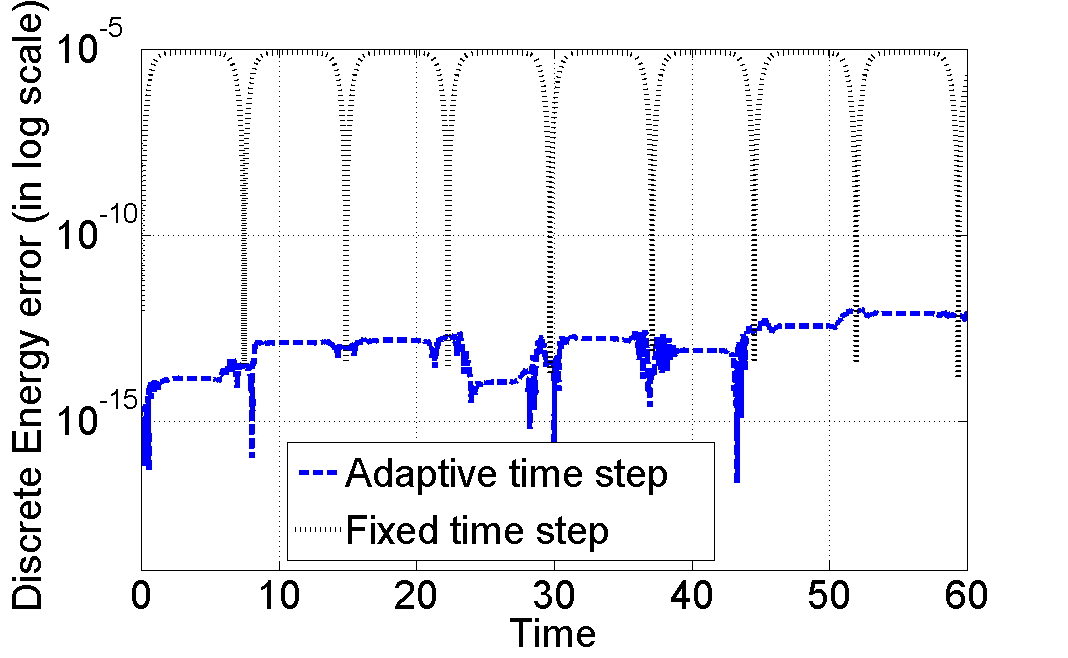}
        \caption{$q(0)=0.995, \ \dot{q}(0)=0$}
        \label{fig:tiger1}
    \end{subfigure}
    \caption{Energy error plots for both fixed time step and adaptive time step algorithms are compared for two different initial conditions. Each figure shows the superior energy performance of the adaptive time step algorithm.}\label{fig:animals1}
\end{figure}
 \\ \\
\textbf{Remark 3.} In \cite{Kane1999Symplectic-energy-momentumIntegrators}, an alternative optimization method has been implemented  where instead of solving the nonlinear coupled equations \eqref{eq:tmpc} and \eqref{eq:tmec}, the following quantity is minimized
\begin{equation}
 [ F(q_k,p_k,h_k,v_k) ]^2 + [ G(q_k,E_k,h_k,v_k) ]^2
\end{equation}
over the variables $v_k$ and $h_k$ with the restriction $h_k>0$. The drawback of using this approach is that numerically it violates the energy evolution equation and the underlying structure is no longer preserved. In fact, due to this optimization approach, the energy plots given in \cite{Kane1999Symplectic-energy-momentumIntegrators} do not clearly convey the advantage of energy-preserving variational integrators.
\begin{figure}[h]
    \centering
    \begin{subfigure}[b]{0.45\textwidth}
        \includegraphics[width=\textwidth]{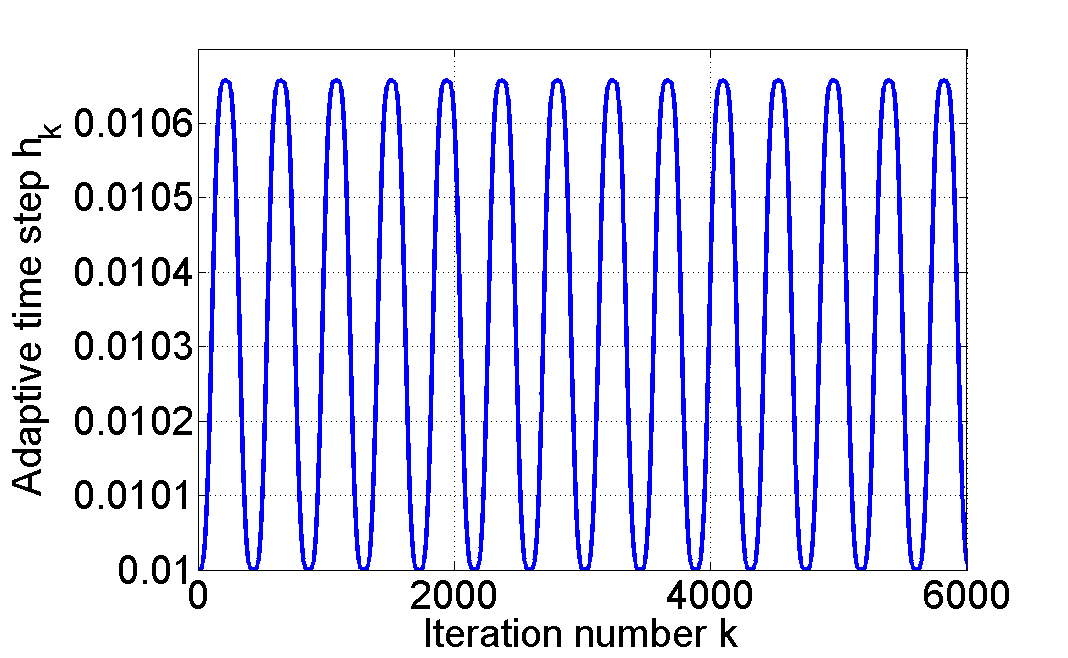}
        \caption{$q(0)=0.74, \ \dot{q}(0)=0$}
        \label{fig:gull2}
    \end{subfigure}
    \begin{subfigure}[b]{0.45\textwidth}
        \includegraphics[width=\textwidth]{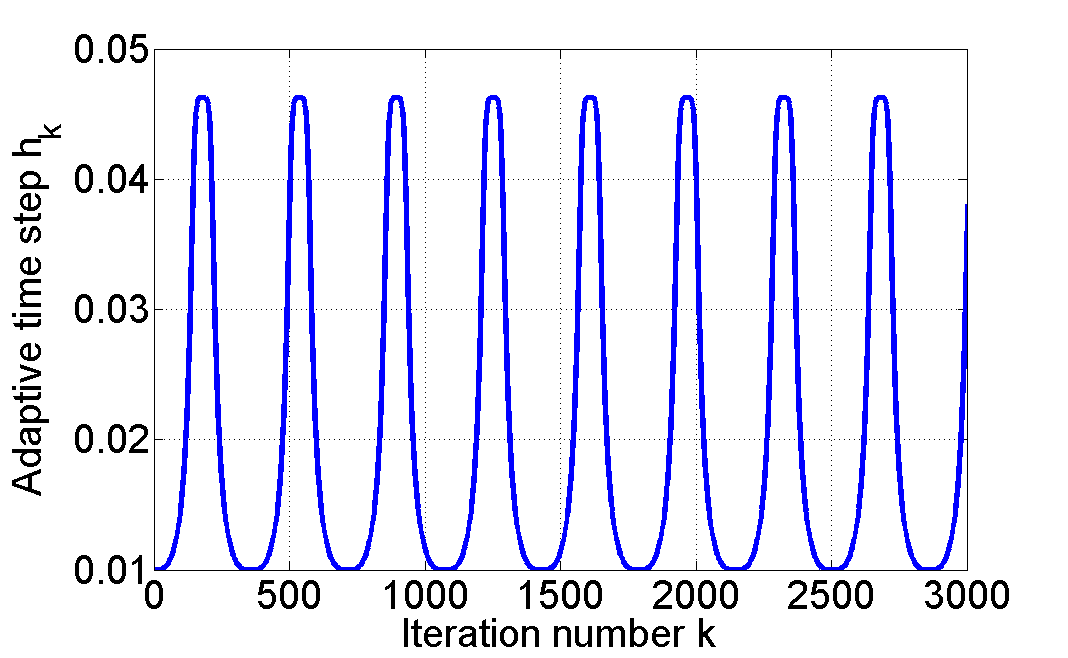}
        \caption{$q(0)=0.995, \ \dot{q}(0)=0$}
        \label{fig:tiger2}
    \end{subfigure}
    \caption{Adaptive time step versus iteration number for both initial conditions.}\label{fig:animals2}
\end{figure} \\
\subsubsection{Initial condition}
We have considered two regions of phase space, similar to the numerical example in \cite{Kane1999Symplectic-energy-momentumIntegrators}, to understand numerical properties of the adaptive time step variational integrator for conservative systems. Since our aim is to use these discrete trajectories as numerical integrators for continuous dynamical systems, instead of starting with two discrete points we consider a continuous-time dynamical system with given initial position $q(0)$ and initial velocity $\dot{q}(0)$ and use the benchmark solution to obtain initial conditions for the adaptive time step variational integrator.
\par 
For a given set of initial conditions, i.e. $(q(0),\dot{q}(0))$, we  first decide the initial time step $h_0$ and then use the benchmark solution to compute discrete configuration $q_1$ at time $t_1=h_0$, configuration at first time step. Thus, we have obtained two discrete points in the extended state space $(t_0,q_0)$ and $(t_1,q_1)$ and using these two discrete points we can find discrete momentum $p_1$ and discrete energy $E_1$. After obtaining $(t_1,q_1,p_1,E_1)$, we can solve the time-marching equations \eqref{eq:pimpf}-\eqref{eq:eexpf} to numerically simulate the dynamical system. 
\subsubsection{Results}
The discrete trajectories for both fixed and adaptive time step algorithms are compared with the benchmark solution in Figure \ref{fig:animals0}. The position $q=\frac{q_k + q_{k+1} }{2}$ and velocity $\dot{q}=\frac{q_{k+1}-q_k}{h_k}$ are computed from the discrete trajectories for both fixed and adaptive time step algorithms and compared with continuous time $q$ and $\dot{q}$. For both initial conditions, discrete trajectories from the adaptive time step and fixed time step  match the benchmark trajectory. 
\par
The energy error plots for both cases show the superior energy behavior of adaptive time step variational integrators for conservative dynamical systems. Instead of using the optimization approach discussed in Remark 3, we have obtained the discrete trajectories by solving the nonlinear coupled equations exactly to preserve the underlying structure. The energy error plots given in Figure \ref{fig:animals1} quantify the difference in energy accuracy for fixed time step and adaptive time step method clearly. The energy-preserving performance was not evident in similar results given in \cite{Kane1999Symplectic-energy-momentumIntegrators} because of the optimization approach used to obtain discrete trajectories instead of solving the implicit equations directly. The energy error comparison in Figure \ref{fig:gull1} shows that the adaptive time step method has energy error  magnitude around $10^{-14}$ whereas the fixed time step method has energy error around $10^{-8}$. In Figure \ref{fig:tiger1} the energy error for fixed time step increases to $10^{-6}$ while the adaptive time step method shows nearly exact energy preservation. Although the magnitude of energy error for fixed time step method is bounded, the magnitude of energy error oscillations depends on where the trajectory lies in the phase space. Thus, for areas in phase space where the magnitude of energy error oscillations is substantial for fixed time step method, the adaptive time step method can be used to preserve the energy of the system more accurately.
\par 
The energy error plots for the adaptive time step algorithm exhibits small jumps in energy error which, we believe, is due to the ill-conditioned nature of the coupled implicit nonlinear equations.
Since the governing implicit equations aren't solved exactly, small numerical errors are introduced at every adaptive time step. These small numerical errors lead to jumps in discrete energy error because of the ill-conditioned nature of the implicit equations.\par
Figure \ref{fig:animals2} shows how the adaptive time step oscillates for both cases. The adaptive time step doesn't increase substantially  compared to the initial time step of $h_0=0.01$ for the first case in Figure \ref{fig:gull2}, while Figure \ref{fig:tiger2} indicates the adaptive time step increases by 4 times the initial time step for the second case. The amplitude of adaptive time step oscillations depends on the region of phase space in which the discrete trajectory lies. The adaptive time step algorithm computes the adaptive time step such that the discrete energy is conserved exactly. There is no upper bound on the size of the adaptive time step, but very large adaptive time step values make the discretization assumption  made in \eqref{eq:adap} erroneous leading to inaccurate discrete trajectories.\\ \\
\textbf{Remark 4.} It is important to understand that adaptive time step variational integrators are fundamentally different from traditional adaptive time-stepping numerical methods which compute the adaptive time step size based on some error criteria. Adaptive time step variational integrators treat time as a discrete dynamic variable and the adaptive time step is computed by solving the extended discrete Euler-Lagrange equations. Thus, the adaptive time step is coupled with the dynamics of the system whereas, for most of the the adaptive time-stepping numerical methods, the step size computation and dynamics of the system are independent of each other.
\begin{figure}[h]
    \centering
    \begin{subfigure}[h]{0.45\textwidth}
        \includegraphics[width=\textwidth]{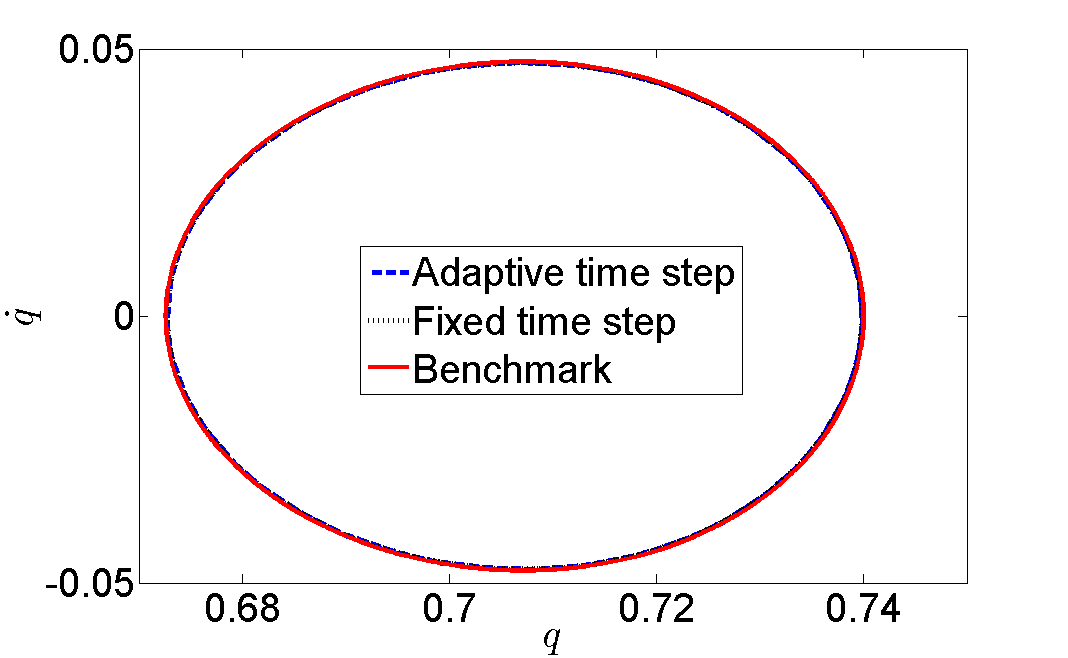}
       \caption{$q(0)=0.74, \ \dot{q}(0)=0$}
        \label{fig:gull4}
    \end{subfigure}
    \begin{subfigure}[h]{0.45\textwidth}
        \includegraphics[width=\textwidth]{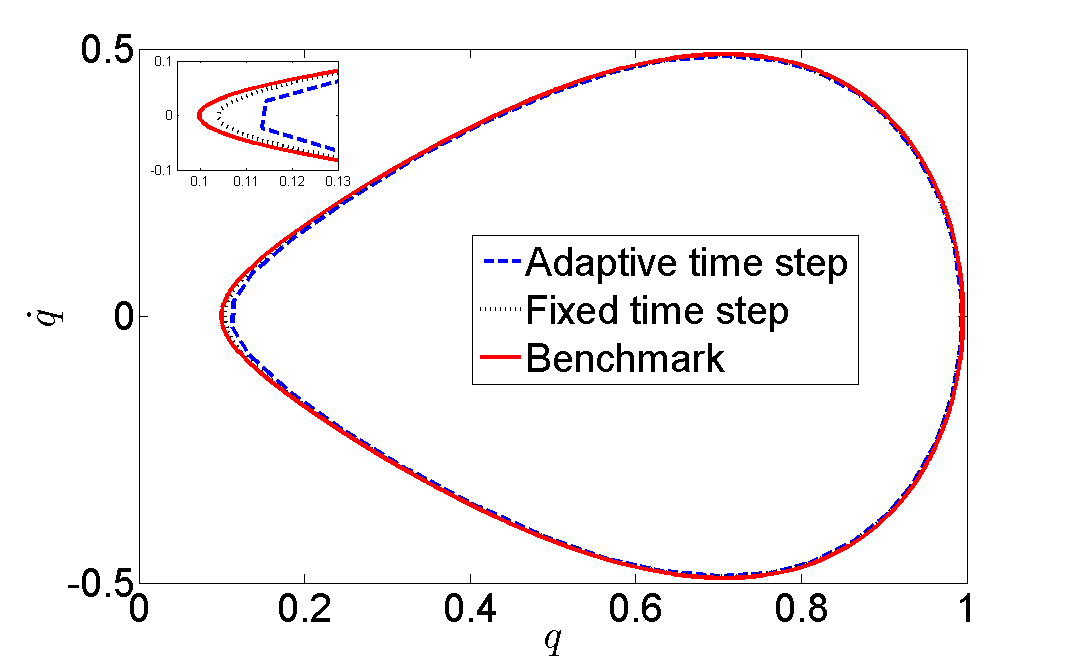}
        \caption{$q(0)=0.995, \ \dot{q}(0)=0$}
        \label{fig:tiger4}
    \end{subfigure}
    \caption{In these plots, an initial time step of $h_0=0.1$ is used  to study the effect of initial time step on the accuracy of discrete trajectories.}\label{fig:animals4}
\end{figure} 
\subsubsection{Effect of initial time step}
From the discrete energy definition it is clear that the initial time step value plays an important role in the adaptive time step algorithm. We study the effect of initial time step on the phase space and energy error plots by simulating the two cases considered in the previous subsection but with a larger initial time step $h_0=0.1$.\par
\begin{figure}[h]
    \centering
    \begin{subfigure}[h]{0.45\textwidth}
        \includegraphics[width=\textwidth]{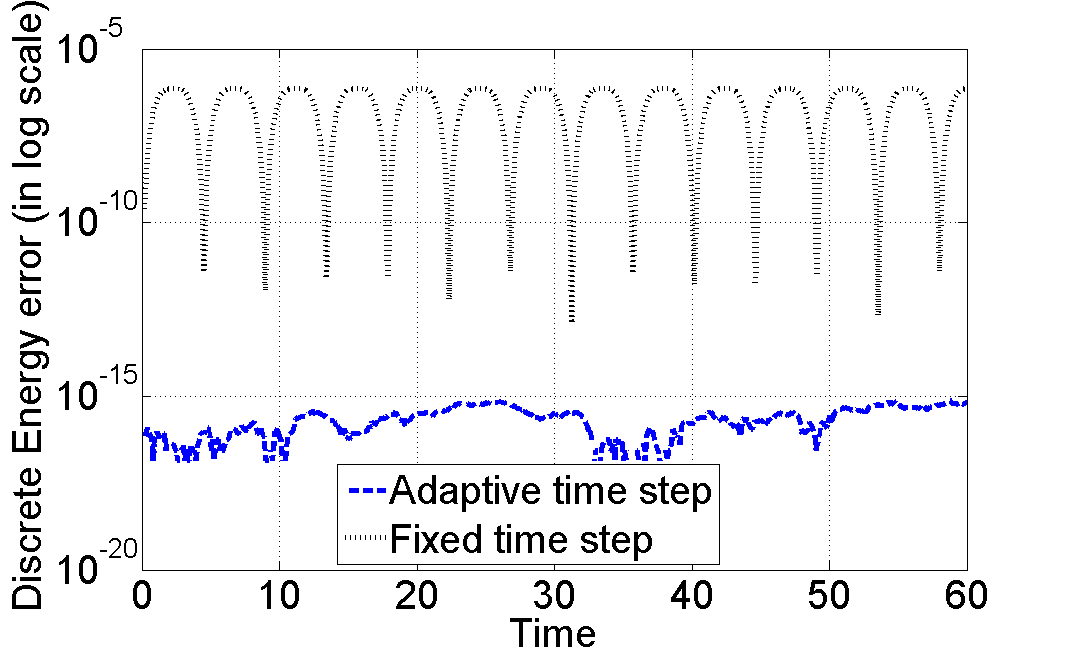}
       \caption{$q(0)=0.74, \ \dot{q}(0)=0$}
        \label{fig:gull5}
    \end{subfigure}
    \begin{subfigure}[h]{0.45\textwidth}
        \includegraphics[width=\textwidth]{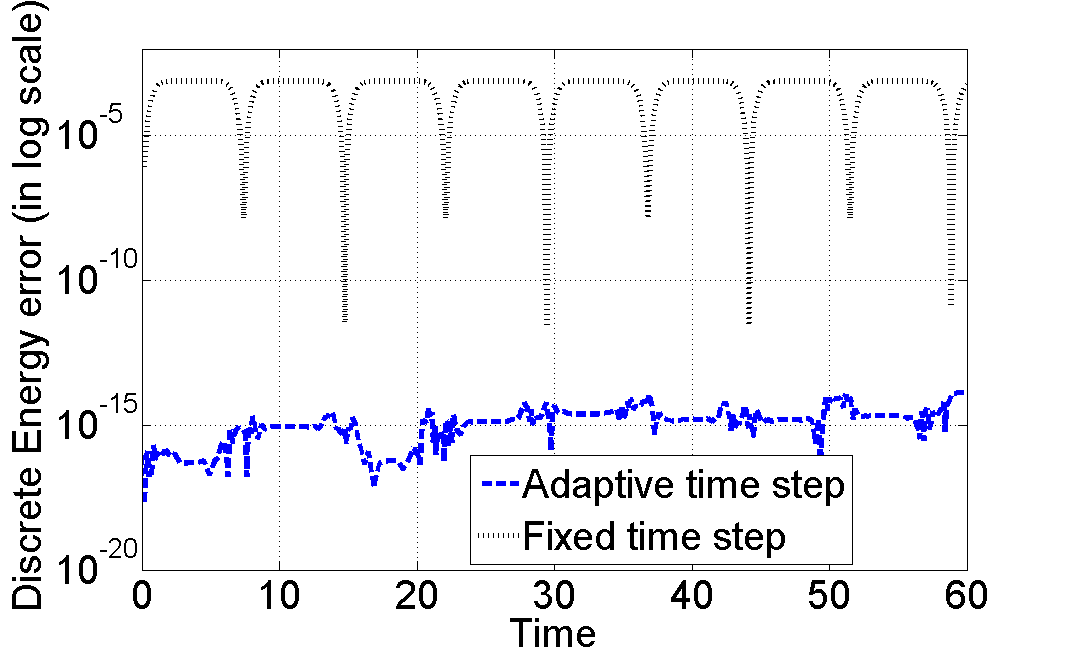}
        \caption{$q(0)=0.995, \ \dot{q}(0)=0$}
        \label{fig:tiger5}
    \end{subfigure}
    \caption{The energy error plots with an initial time step of $h_0=0.1$.}\label{fig:animals5}
\end{figure}
\begin{figure}[h]
    \centering
    \begin{subfigure}[h]{0.45\textwidth}
        \includegraphics[width=\textwidth]{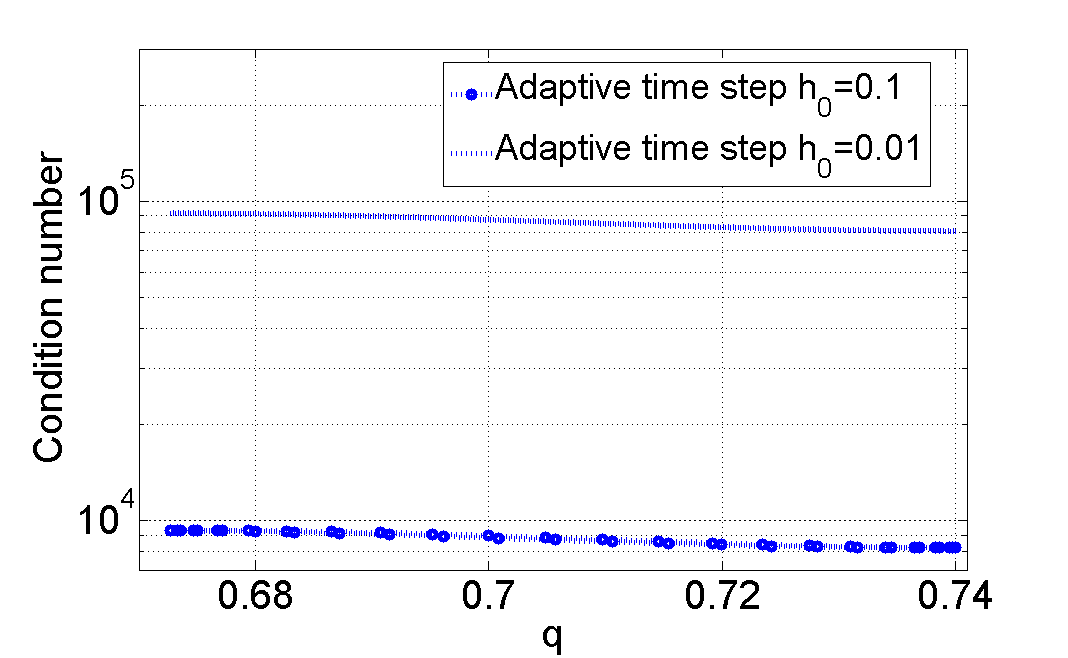}
       \caption{$q(0)=0.74, \ \dot{q}(0)=0$}
        \label{fig:gull6}
    \end{subfigure}
    \begin{subfigure}[h]{0.45\textwidth}
        \includegraphics[width=\textwidth]{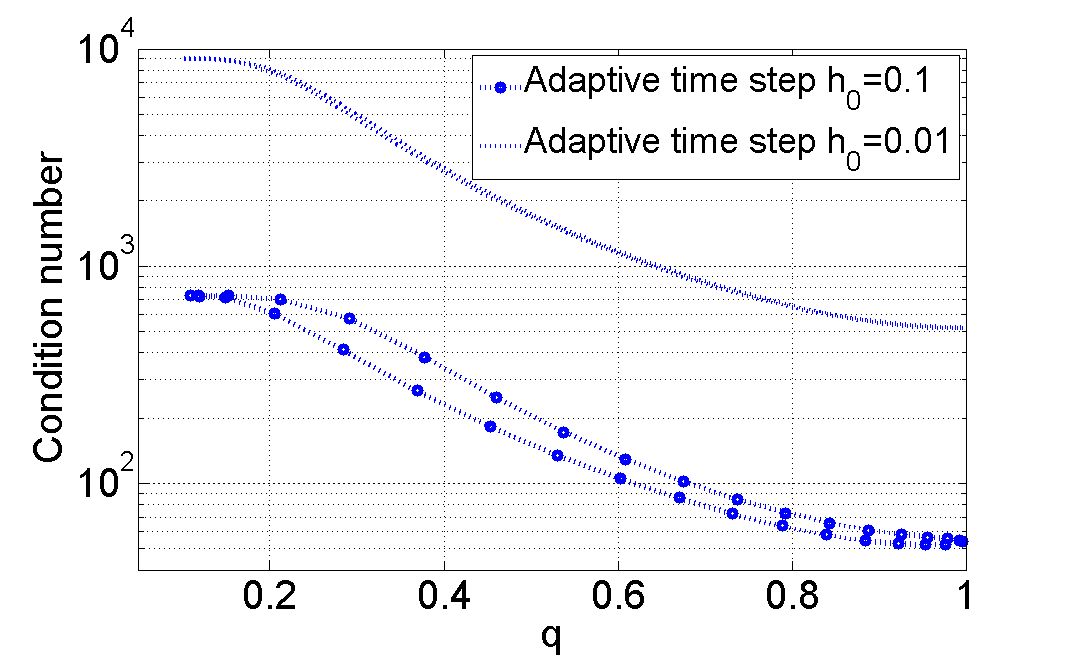}
        \caption{$q(0)=0.995, \ \dot{q}(0)=0$}
        \label{fig:tiger6}
    \end{subfigure}
    \caption{The condition numbers increase with decrease in initial time step in the case of the adaptive time step algorithm.}\label{fig:animals6}
\end{figure}

The phase space trajectories shown in Figure \ref{fig:animals4} show that even with an initial time step of $h_0=0.1$ discrete trajectories from both fixed and adaptive time step show good agreement with the benchmark solution. In Figure \ref{fig:gull4}, the discrete trajectories lie on top of the benchmark solution for the first set of initial conditions. In Figure \ref{fig:tiger4}, the fixed and adaptive time step discrete trajectories give slightly inaccurate results near the turning point. The discrete energy error plots in Figure \ref{fig:animals5} show that for fixed time step variational integrators, the discrete energy errors increase with increase in the time step size but for the adaptive time step variational integrators, increasing the time step size leads to more accurate discrete energy behavior. This unexpected behavior is due to the ill-conditioned nature of the implicit extended discrete Euler-Lagrange equations, which become more ill-conditioned for smaller time steps. The plots of condition number in Figure \ref{fig:animals6} show that the implicit equations become more ill-conditioned as the initial time step value is decreased. \par
It is important to note that for a conservative system, the continuous-time trajectory preserves the continuous energy which is different from the discrete energy that adaptive time step variational integrators are constructed to preserve. This explains why, despite the superior energy behavior in Figure \ref{fig:tiger5} compared to Figure \ref{fig:tiger1}, the discrete trajectory in Figure \ref{fig:tiger4} is less accurate than the discrete trajectory in Figure \ref{fig:tiger0}. We know that as the time step value tends to zero the discrete energy and continuous energy become equal but the condition number analysis and the energy error plots reveal that smaller initial time steps for adaptive time variational integrators lead to an increase in energy error. \emph{Thus, there is a trade-off between preserving discrete energy and ensuring accuracy when choosing an initial time step for the adaptive time step variational integrators.}
 \\ 
\textbf{Remark 5.} The discrete energy error plotted in Figure \ref{fig:animals1} and Figure \ref{fig:animals5} is different from our traditional idea of energy error. We usually define energy error as the difference between the energy of the continuous-time system and the energy obtained from the discrete trajectories. This traditional energy error can be broken down into \emph{discrete energy error} and \emph{discretization error}. The discrete energy error is the error in preserving the discrete energy of the extended discrete Lagrangian system. The discretization error is the error incurred by discretizing a continuous-time system. Thus, the discretization error is the difference between the continuous energy and the discrete energy that our integrators aim to preserve, whereas the discrete energy error is the error between the true and computed discrete energy.
\\ 
\textbf{Remark 6.} For a conservative system, we expect discrete energy to be constant and thus the discretization error is also constant. Since this constant discretization error is orders of magnitude larger than the discrete energy error, traditional energy error plots do not show the advantages of using adaptive time-stepping. We evaluate the performance of variational integrators by comparing how well these integrators preserve the discrete energy.
\subsection{Dissipative Example}
\label{sc:52}
We consider a damped harmonic oscillator in order to better understand the numerical behavior of the adaptive time step variational integrator for forced Lagrangian systems. The (continuous) Lagrangian for the single degree of freedom system is 
\begin{equation}
L (q, \dot{q}) = \frac{1}{2}m \dot{q}^2 - \frac{1}{2}kq^2
\end{equation}
and the dissipative force is
\begin{equation}
f= -c\dot{q}
\end{equation}
where $m$ is the mass, $k$ is the stiffness and $c$ is the damping parameter of the single degree of freedom system. For the discrete Lagrangian $L_d$, we use the midpoint rule which gives
\begin{equation}
L_d (t_k,q_k,t_{k+1},q_{k+1}) = (t_{k+1}-t_k) \  L\left (\frac{q_k + q_{k+1}}{2},\frac{q_{k+1}-q_k}{t_{k+1}-t_k}\right)
\end{equation}
Similarly, we can write the discrete force $f_d^{\pm}$ as
\begin{equation}
f_d^{\pm} = -\frac{1}{2}c(t_{k+1}-t_k)  \left( \frac{q_{k+1}-q_k}{t_{k+1}-t_k} \right)
\end{equation}
and the corresponding power term $g_d^{\pm}$ is 
\begin{equation}
g_d^{\pm}=-f_d^{\pm} \left( \frac{q_{k+1}-q_k}{t_{k+1}-t_k} \right) = \frac{1}{2}c(t_{k+1}-t_k)  \left( \frac{q_{k+1}-q_k}{t_{k+1}-t_k} \right)^2
\end{equation}
The discrete momentum $p_{k+1}$ and discrete energy $E_{k+1}$ expressions are  
\begin{equation}
\begin{aligned}
p_{k+1}&=D_4L_d(t_k,q_k,t_{k+1},q_{k+1}) +  f_d^+ \\ &= m\left ( \frac{q_{k+1}-q_k}{t_{k+1}-t_k} \right)- k(t_{k+1}-t_k) \left( \frac{q_k + q_{k+1}}{4} \right) - c\left( \frac{q_{k+1}-q_k}{2}\right)
\end{aligned}
\label{eq:pk1_nex}
\end{equation}

\begin{equation}
\begin{aligned}
E_{k+1} &= - D_3L_d(t_{k},q_{k},t_{k+1},q_{k+1}) -  g_d^+ \\ &= \frac{1}{2}m \left( \frac{q_{k+1}-q_k}{t_{k+1}-t_k} \right)^2 + \frac{1}{2}k \left( \frac{q_k + q_{k+1}}{2} \right)^2  - c \left( \frac{(q_{k+1}-q_k)^2}{t_{k+1}-t_k} \right)
\label{eq:ek1_nex}
\end{aligned}
\end{equation}
\begin{figure}[h]
    \centering
    \begin{subfigure}[h]{0.45\textwidth}
        \includegraphics[width=\textwidth]{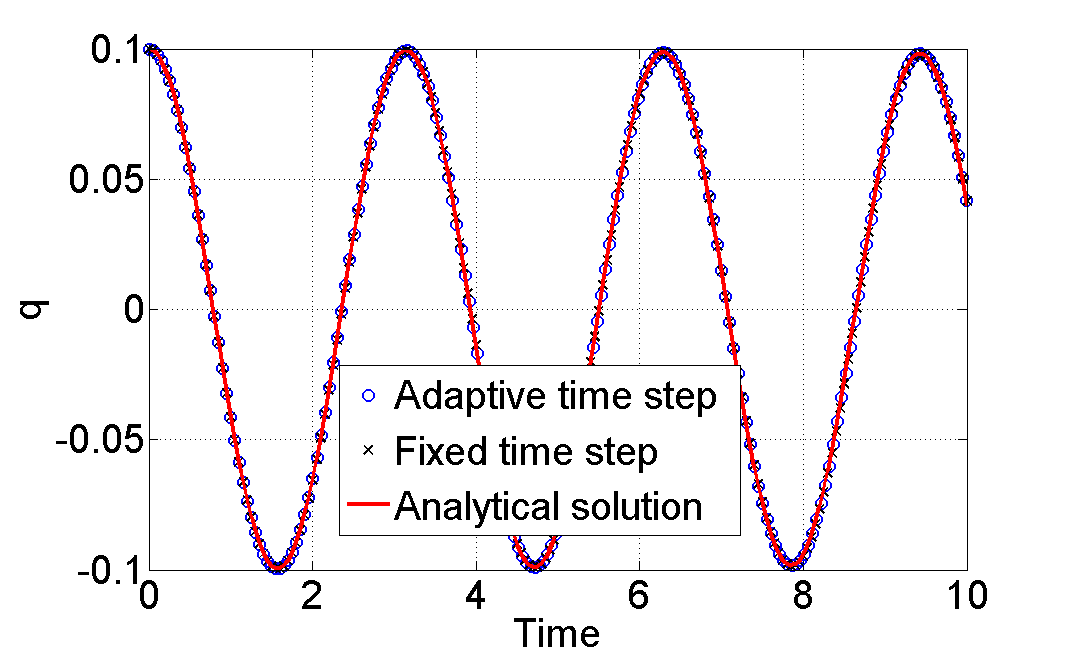}
        \caption{$\zeta=0.001$}
        \label{fig:gull7}
    \end{subfigure}
    \begin{subfigure}[h]{0.45\textwidth}
        \includegraphics[width=\textwidth]{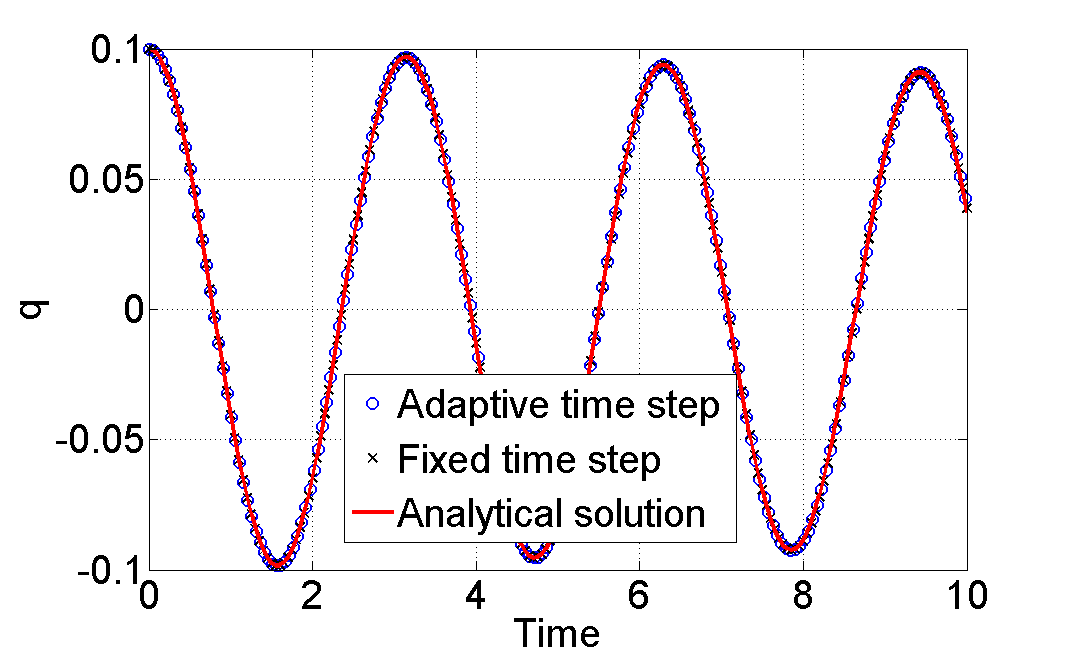}
        \caption{$\zeta=0.005$}
        \label{fig:tiger7}
    \end{subfigure}
    \begin{subfigure}[h]{0.45\textwidth}
        \includegraphics[width=\textwidth]{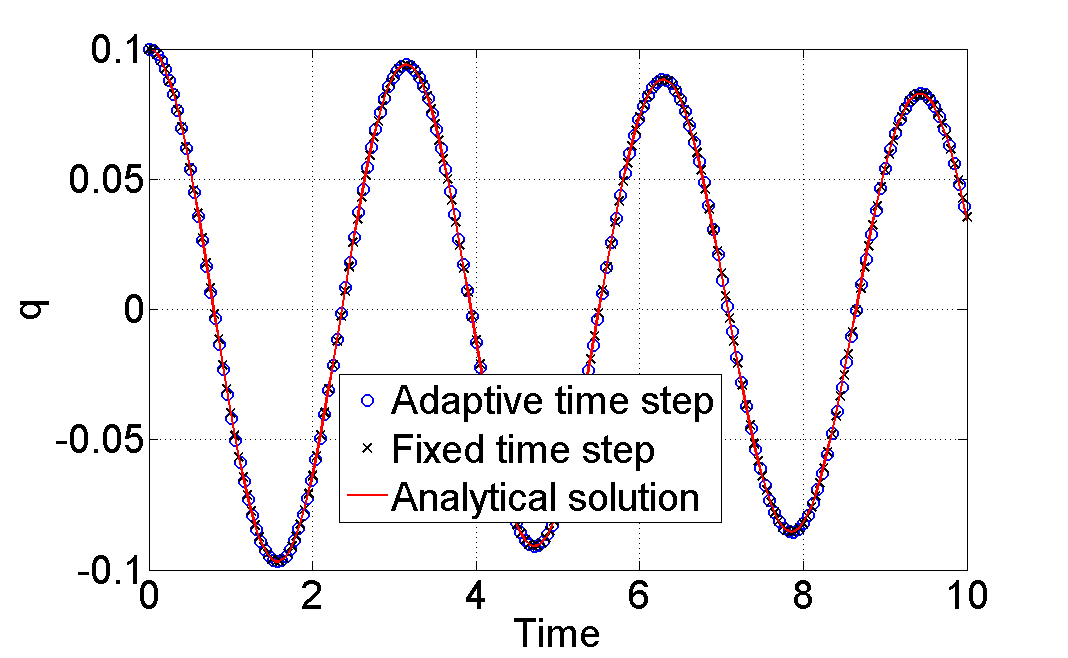}
        \caption{$\zeta=0.01$}
        \label{fig:mouse7}
    \end{subfigure}
    \caption{Three damping ratio values are studied for the spring mass damper system. Discrete trajectories for both fixed time step and adaptive time step variational integrators are plotted and compared with the analytical solution. The analytical solution is used to prescribe initial conditions for an initial time step $h_0=0.01$ and natural frequency $\omega_n=2 \  rad/s$. }\label{fig:animals7}
\end{figure} \\
For given $(t_k,q_k,E_k,p_k)$, the time-marching implicit equations are obtained by substituting the discrete Lagrangian and discrete force expressions into \eqref{eq:pimpf} and \eqref{eq:eimpf}
\begin{equation}
m\left ( \frac{q_{k+1}+q_k}{t_{k+1}-t_k} \right)  + k(t_{k+1}-t_k)\left ( \frac{q_{k+1}+q_k}{4} \right) - c(t_{k+1}-t_k)\left(\frac{q_{k+1}-q_k}{2}\right) =p_k
\end{equation}
\begin{equation}
\frac{1}{2}m\left ( \frac{q_{k+1}-q_k}{t_{k+1}-t_k} \right)^2 + \frac{1}{2}c(t_{k+1}-t_k)\left ( \frac{q_{k+1}-q_k}{t_{k+1}-t_k} \right)^2 + \frac{1}{2}k \left ( \frac{q_{k+1}+q_k}{2} \right)^2 = E_k
\end{equation}
The above two coupled nonlinear equations in $q_{k+1}$ and $t_{k+1}$ are solved with the restriction \  \ \  \ \ \  $t_{k+1}>t_k$ and substituted in \eqref{eq:pk1_nex} and \eqref{eq:ek1_nex} to obtain the discrete momentum $p_{k+1}$ and discrete energy $E_{k+1}$ for the next step. We re-write the above time-marching equations in terms of $h_k=t_{k+1}-t_k$ and $v_k=\left( \frac{q_{k+1}-q_k}{t_{k+1}-t_k} \right)$
\begin{equation}
F(q_k,p_k,h_k,v_k) = mv_k + \frac{h_k}{4} ( 2q_k + h_kv_k ) + \frac{1}{2}ch_kv_k - p_k = 0 
\label{eq:tmp}
\end{equation}
\begin{equation}
G(q_k,E_k,h_k,v_k) = \frac{1}{2}mv_k^2 + \frac{1}{2}ch_kv_k^2 + \frac{1}{2}k \left( q_k + \frac{h_kv_k}{2} \right)^2 - E_k = 0 
\label{eq:tmE}
\end{equation}
\textbf{Remark 7.} Since the Lagrangian and the forcing for this example are both time-independent, we can replace $t_{k+1}-t_k$ by the $k^{th}$ adaptive time step $h_k$ as shown above. For time-dependent mechanical systems with either time-dependent Lagrangian or time-dependent forcing, this simplification cannot be made and the implicit equations must be solved for $t_{k+1}$.
\begin{figure}[h]
    \centering
    \begin{subfigure}[h]{0.45\textwidth}
        \includegraphics[width=\textwidth]{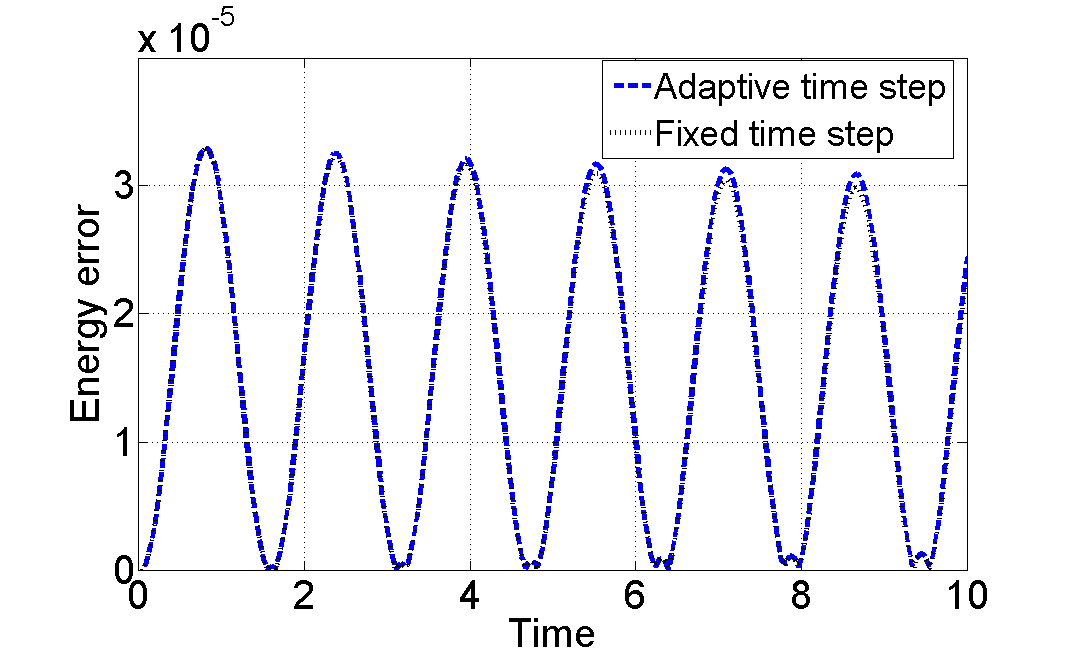}
        \caption{$\zeta=0.001$}
        \label{fig:gull8}
    \end{subfigure}
    \begin{subfigure}[h]{0.45\textwidth}
        \includegraphics[width=\textwidth]{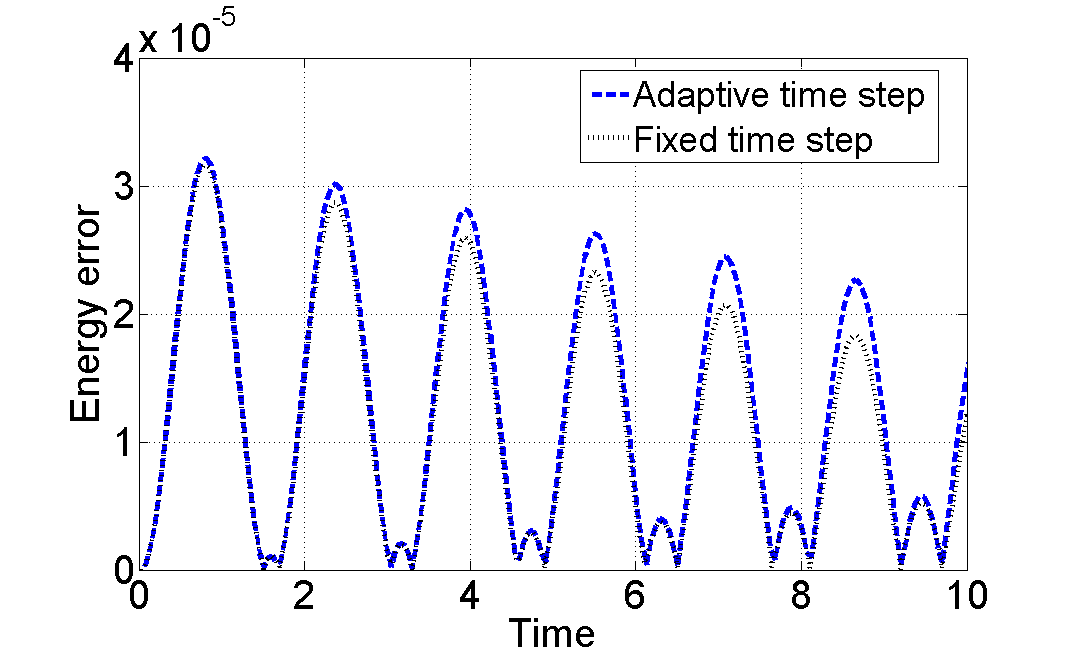}
        \caption{$\zeta=0.005$}
        \label{fig:tiger8}
    \end{subfigure}
    \begin{subfigure}[h]{0.45\textwidth}
        \includegraphics[width=\textwidth]{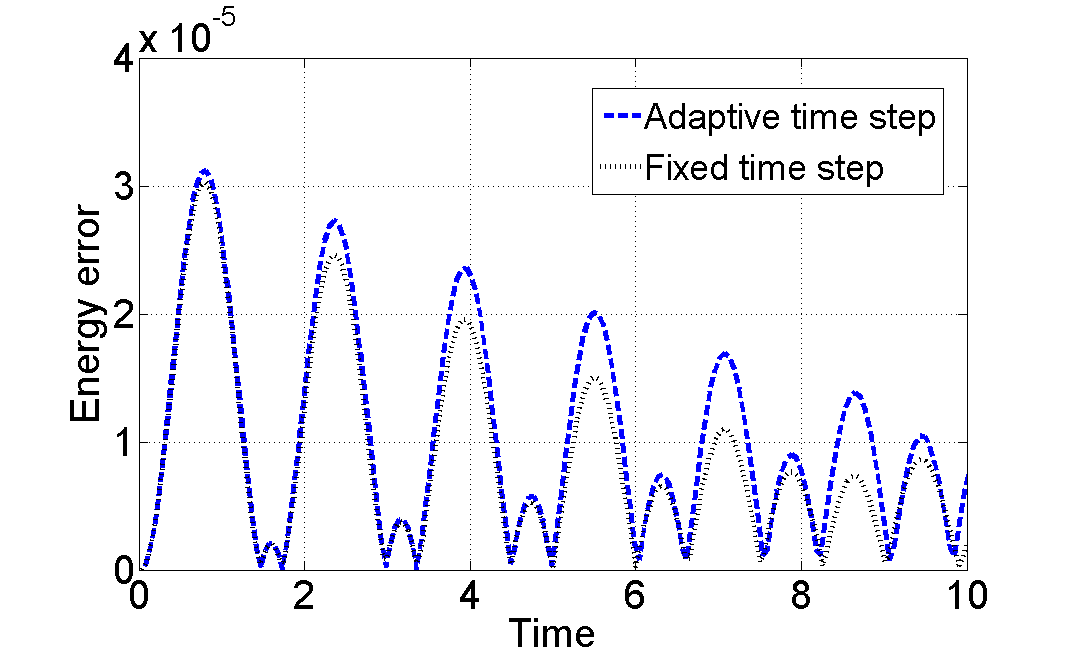}
        \caption{$\zeta=0.01$}
        \label{fig:mouse8}
    \end{subfigure}
    \caption{Energy error for fixed time step and adaptive time step variational integrators are compared for three cases. Analytical solution at the discrete time instant is used to compute the continuous energy. }\label{fig:animals8}
\end{figure}
\begin{figure}[h]
    \centering
    \includegraphics[width=\textwidth]{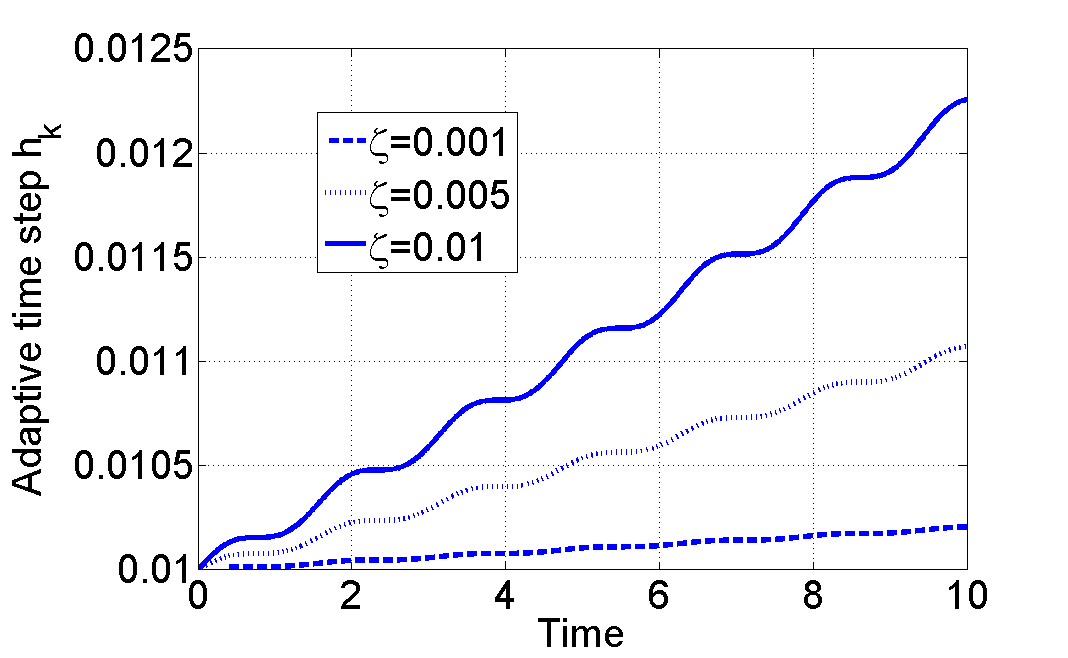}
    \caption{Adaptive time step  versus the time for the energy-preserving variational integrator.}
    \label{fig:mesh1}
\end{figure}  
\subsubsection{Results}
We have studied the damped simple harmonic oscillator for three small damping values of the damping ratio $\zeta = \frac{c}{2\sqrt{km}}$ for a single natural frequency $\omega_n=\sqrt{\frac{k}{m}}=2 \  rad/s$ to understand the numerical properties of adaptive time step variational integrators derived for forced systems in Section \ref{s:4}. Just like the conservative case, our aim is to simulate the continuous-time dynamical system using discrete trajectories obtained from the adaptive time step variational integrator. The discrete  trajectories from both fixed and adaptive time step algorithms are compared in Figure \ref{fig:animals7}. Both are nearly indistinguishable from the analytical solution for all three cases. \par
%
%
%
%
The energy error plots in Figure \ref{fig:animals8} show how both adaptive and fixed time step variational integrators  start with same energy accuracy for all three cases but, as we march forward in time, the fixed time step variational integrator outperforms the adaptive time step variational integrator. \emph{The amplitude of the energy error oscillations for the fixed time step algorithm decreases faster than it does for the adaptive time step algorithm which suggests that for long-time simulations the energy behavior of fixed time step variational integrator is better than the adaptive time step variational integrator.} This is contrary to what we expected because the adaptive step variational integrator solves an additional discrete energy evolution equation to capture the change in energy of the forced system accurately.\par
%
%
%
These unexpected results can be understood by looking at the two components of the energy error discussed in the Remark 5. Due to exact preservation of discrete energy, the discrete energy error for adaptive time step variational integrators is orders of magnitude lower than it is for the fixed time step variational integrator. Unlike the  conservative system example considered \ \  in Section \ref{sc:51}, the continuous energy and the corresponding discrete energy, for this dissipative system, are not constant. Thus, the energy errors are computed by comparing the continuous energy with the discrete counterpart. The discrete energy for forced Lagrangian systems has terms accounting for virtual work done by the external force during the adaptive time step and  hence the adaptive time step variational integrators are preserving a discrete quantity which is not analogous to the continuous time energy. Since the difference between continuous and discrete energy is orders of magnitude larger than the  discrete energy error of  the variational integrator, the resulting energy error plots do not reflect the advantage of using adaptive variational integrators over fixed time step variational integrators. \par
Another reason behind the higher energy error for adaptive time step variational integrators is the monotonically increasing adaptive time step shown in Figure \ref{fig:mesh1}. The velocity approximation $\dot{q}\approx \frac{q_{k+1}-q_k}{h_k}$ used in computing the discrete energy becomes more inaccurate as the adaptive time step increases. As we go forward in time, the adaptive time step $h_k$ keeps on increasing leading to higher energy error for adaptive time step variational integrators. Thus, the magnitude of energy error for adaptive time does not decrease as quickly as it does for fixed time step variational integrators.\par
In Figure \ref{fig:mesh1} the adaptive time step evolution over time for all three damping parameter values is plotted. For all three cases, the adaptive time step was found to be monotonically increasing. This is not good for a numerical algorithm as eventually it would lead to numerical instability. We have also studied the damped harmonic oscillator system for negative damping parameter values and the results for those systems showed a uniformly decreasing adaptive time step. Thus, there seems to be some inverse relation between the rate of change of energy and the rate of change of the adaptive time step.
\section{Conclusions and Future work}
\label{s:6}
In this work we have presented adaptive time step variational integrators for time-dependent mechanical systems with forcing. We have incorporated forcing into the extended discrete mechanics framework so that the resulting discrete trajectories can be used as numerical integrators for Lagrangian systems with forcing. The paper first presented the Lagrange-d'Alembert principle in the extended Lagrangian mechanics framework and then derived the extended forced discrete Euler-Lagrange equations from the discrete Lagrange-d'Alembert principle. We demonstrated a general method to construct adaptive time step variational integrators for Lagrangian systems with forcing through a damped harmonic oscillator example. The results from the numerical example showed that the adaptive time step algorithm works well for cases where the rate of change of energy of the system is very slow. The adaptive time step for the dissipative system was found to be monotonically increasing which makes the algorithm unsuitable for long-time simulation. 
\par
We have also presented results for  a nonlinear conservative system by solving the discrete equations exactly, as opposed to the optimization approach suggested in \cite{Kane1999Symplectic-energy-momentumIntegrators}. The energy error results show the advantage of solving discrete equations exactly for adaptive time step variational integrators. We have studied the effect of initial time step on energy error and phase space trajectories and also shown how the discrete equations become more ill-conditioned as the initial time step becomes smaller.\par 
In future work, we would like to study the connection between the rate of change of energy and the size of the adaptive time step which is evident in the example of the damped harmonic oscillator. It would also be desirable to investigate the numerical performance of variational integrators for time-dependent Lagrangian systems. \\ 
\\
\textbf{Declaration of interest:} None.
\\ 
\\
\textbf{Funding:} This research did not receive any specific grant from funding agencies in the public, commercial, or
not-for-profit sectors.

\bibliographystyle{vancouver}
\bibliography{main}
\end{document}